\documentclass[graybox]{svmult}
\usepackage{amssymb, amsfonts, amsbsy, latexsym, epsfig,
  color,amsmath,pstool}

\usepackage{mathptmx}       
\usepackage{helvet}         
\usepackage{courier}        
\usepackage{type1cm}        
\usepackage{makeidx}         
\usepackage{graphicx}        
\usepackage{multicol}        
\usepackage[bottom]{footmisc}

\usepackage{srcltx}
\usepackage{enumerate}
\usepackage[T1]{fontenc}
\usepackage[latin1]{inputenc}
\usepackage{graphicx}
\usepackage{psfrag}
\usepackage{pst-all}
\usepackage{pstricks-add}
\usepackage[caption=false]{subfig}
\usepackage{cite}
\usepackage{url}

\usepackage{mathrsfs}

\allowdisplaybreaks[3]

\makeindex             

\DeclareMathOperator{\supp}{supp} %
\DeclareMathOperator{\diagonal}{diag} %
\DeclareMathOperator*{\esssup}{ess\,sup} %
\DeclareMathOperator*{\essinf}{ess\,inf} %

\newcommand{\diag}[1]{\diagonal{(#1)}}
\newcommand{\sign}[1]{\mathrm{sgn}(#1)}

\newcommand*{\numbersys}[1]{\ensuremath{\mathbb{#1}}}

\newcommand*{\R}{\numbersys{R}}

\newcommand*{\Z}{\numbersys{Z}}

\newcommand*{\N}{\numbersys{N}}

\newcommand{\itvcc}[2]{\ensuremath{\left[{#1},{#2}\right]}} %
\newcommand{\itvco}[2]{\ensuremath{\left[{#1},{#2}\right)}} %
\newcommand{\abs}[1]{\ensuremath{\left\lvert#1\right\rvert}}
\newcommand{\abssmall}[1]{\ensuremath{\lvert#1\rvert}}

\newcommand{\norm}[2][]{\ensuremath{\left\lVert#2\right\rVert_{#1}}}

\newcommand{\innerprod}[3][]{\ensuremath{\left\langle #2,#3\right\rangle_{\! #1}}}

\newcommand{\setprop}[2]{\ensuremath{\left\lbrace{#1} : {#2}\right\rbrace}}

\usepackage{xspace}
\newcommand{\ie}{i.e.,\xspace} 
\newcommand{\eg}{e.g.,\xspace} 

\newcommand{\scfac}{2^j}

\newcommand{\ceilsmall}[1]{\lceil #1 \rceil}

\newcommand\cC{{\mathcal{C}}}

\newcommand\cE{{\mathcal{E}}}

\newcommand\cP{{\mathcal{P}}}
\newcommand\cR{{\mathcal{R}}}

\newcommand\SH{SH}

\usepackage[colorlinks=true, linkcolor=black, citecolor=black,
filecolor=black, pagecolor=black, urlcolor=black, bookmarks=true,
bookmarksopen=true, bookmarksopenlevel=3, plainpages=false,
pdfpagelabels=true]{hyperref}
\hypersetup{
pdfview={FitH},
pdfstartview={FitH},
 bookmarks={true},
pdfauthor = {Gitta Kutyniok, Jakob Lemvig, and Wang-Q Lim},
pdftitle = {Compactly Supported Shearlets},
pdfsubject = {shearlet theory},
pdfkeywords = {shearlet, wavelet, anisotropy},
pdfcreator = {LaTeX with hyperref package},
pdfproducer = {dvips + ps2pdf}}

\newcommand{\jl}[1]{#1}

\newcommand{\jll}[1]{#1}

\newcommand{\gggk}[1]{#1}
\newcommand{\wql}[1]{#1}
\begin{document}

\title*{Compactly Supported Shearlets}
\author{Gitta Kutyniok, Jakob Lemvig, and Wang-Q Lim}
\institute{Gitta Kutyniok \at Institute of Mathematics, University of
  Osnabr\"uck, 49069 Osnabr\"uck, Germany,
  \email{kutyniok@uni-osnabrueck.de} \and
Jakob Lemvig \at Institute of Mathematics, University of
  Osnabr\"uck, 49069 Osnabr\"uck, Germany,
  \email{jlemvig@uni-osnabrueck.de} \and
Wang-Q Lim \at Institute of Mathematics, University of
  Osnabr\"uck, 49069 Osnabr\"uck, Germany,
  \email{wlim@uni-osnabrueck.de} }


%
%
\maketitle

\abstract{Shearlet theory has become a central tool in analyzing and
  representing 2D data with anisotropic features. Shearlet systems are
  systems of functions generated by one single generator with
  parabolic scaling, shearing, and translation operators applied to
  it, in much the same way wavelet systems are dyadic scalings and
  translations of a single function, but including a precise
  control of directionality. Of the many directional representation
  systems proposed in the last decade, shearlets are among the most
  versatile and successful systems. The reason for this being an
  extensive list of desirable properties: shearlet systems can be
  generated by one function, they provide precise resolution of
  wavefront sets, they allow compactly supported analyzing elements,
  they are associated with fast decomposition algorithms, and they
  provide a unified treatment of the continuum and the digital realm.\\
  \indent The aim of this paper is to introduce some key concepts in
  directional representation systems and to shed some
  light on the success of shearlet systems as directional
  representation systems. In particular, we will give an overview of the
  different paths taken in shearlet theory with focus on separable and
  compactly supported shearlets in 2D and 3D. We will present constructions
  of compactly supported shearlet frames in those dimensions as well as discuss
  recent results on the ability of compactly supported shearlet frames
  satisfying weak decay, smoothness, and directional moment conditions to
  provide optimally sparse approximations of cartoon-like images in 2D \gggk{as well as in 3D}.
  Finally, we will show that these compactly supported shearlet systems provide
  optimally sparse approximations of an even generalized model of cartoon-like
  images comprising of $C^2$ functions that are smooth apart from
  piecewise $C^2$ discontinuity edges.}

\abstract*{Shearlet theory has become a central tool in analyzing and
  representing 2D data with anisotropic features. Shearlet systems are
  systems of functions generated by one single generator with
  parabolic scaling, shearing, and translation operators applied to
  it, in much the same way wavelet systems are dyadic scalings and
  translations of a single function, but including a precise
  control of directionality. Of the many directional representation
  systems proposed in the last decade, shearlets are among the most
  versatile and successful systems. The reason for this being an
  extensive list of desirable properties: shearlet systems can be
  generated by one function, they provide precise resolution of
  wavefront sets, they allow compactly supported analyzing elements,
  they are associated with fast decomposition algorithms, and they
  provide a unified treatment of the continuum and the digital realm.\\
  \indent The aim of this paper is to introduce some key concepts in
  directional representation systems and to shed some
  light on the success of shearlet systems as directional
  representation systems. In particular, we will give an overview of the
  different paths taken in shearlet theory with focus on separable and
  compactly supported shearlets in 2D and 3D. We will present constructions
  of compactly supported shearlet frames in those dimensions as well as discuss
  recent results on the ability of compactly supported shearlet frames
  satisfying weak decay, smoothness, and directional moment conditions to
  provide optimally sparse approximations of cartoon-like images in 2D \gggk{as well as in 3D}.
  Finally, we will show that these compactly supported shearlet systems provide
  optimally sparse approximations of an even generalized model of cartoon-like
  images comprising of $C^2$ functions that are smooth apart from
  piecewise $C^2$ discontinuity edges.}

\section{Introduction}

Recent advances in modern technology have created a brave new world
of enormous, multi-dimensional data structures. In medical imaging,
seismic imaging, astronomical imaging, computer vision, and video
processing, the capabilities of modern computers and high-precision
measuring devices have generated 2D, 3D, and even higher dimensional
data sets of sizes that were infeasible just a few years ago. The
need to efficiently handle such diverse types and huge amounts of
data initiated an intense study in developing efficient multivariate
encoding methodologies in the applied harmonic analysis research
community.

In medical imaging, \eg CT lung scans, the discontinuity curves of the
image are important specific features since one often wants to
distinguish between the image `objects' (\eg the lungs) and the
`background'; that is, it is important to precisely capture the
\emph{edges}. This observation holds for various other applications
than medical imaging and illustrates that important classes of
multivariate problems are governed by \emph{anisotropic features}.
Moreover, in high-dimensional data most information is typically
contained in lower-dimensional embedded manifolds, thereby also
presenting itself as anisotropic features. The anisotropic structures
can be distinguished by location and orientation/direction which
indicates that our way of analyzing and representing the data should
capture not only location, but also directional information.

In applied harmonic analysis, data is typically modeled in a
continuum setting as square-integrable functions or, more generally,
as distributions. Recently, a novel directional representation
system -- so-called shearlets -- has emerged which provides a
unified treatment of such continuum models as well as digital
models, allowing, for instance, a precise resolution of wavefront
sets, optimally sparse representations of cartoon-like images, and
associated fast decomposition algorithms. Shearlet systems are
systems generated by one single generator with parabolic scaling,
shearing, and translation operators applied to it, in the same way
wavelet systems are dyadic scalings and translations of a single
function, but including a directionality characteristic owing to the
additional shearing operation (and the anisotropic scaling).

The aim of this survey paper is to introduce the key concepts in
directional representation systems and, in particular, to shed some
light on the success of shearlet systems. Moreover, we will give an
overview of the different paths taken in shearlet theory with focus
on separable and compactly supported shearlets, since these systems
are most well-suited for applications in, \eg image processing and
the theory of partial differential equations.

\subsection{Directional Representation Systems}
\label{sec:direct-repr-syst}
In recent years numerous approaches for efficiently representing
directional features of two-dimensional data have been proposed. A
perfunctory list includes: \emph{steerable pyramid} by Simoncelli
\textit{et al.} \cite{steer1992}, \emph{directional filter banks} by
Bamberger and Smith \cite{directional1992}, \emph{2D directional
  wavelets} by Antoine \textit{et al.} \cite{twoDwavelet1993},
\emph{curvelets} by Cand\`es and Donoho \cite{CD99},
\emph{contourlets} by Do and Vetterli \cite{DV05}, \emph{bandlets} by
LePennec and Mallat \cite{bandelets}, and \emph{shearlets} by Labate,
Weiss, and two of the authors \cite{LLKW05}. Of these, shearlets are
among the most versatile and successful systems which owes to the many
desirable properties possessed by shearlet systems: they are generated
by one function, they provide optimally sparse approximation of
so-called cartoon-like images, they allow compactly supported
analyzing elements, they are associated with fast decomposition
algorithms, and they provide a unified treatment of continuum and
digital data.

Cartoon-like images are functions that are $C^2$ apart from $C^2$
singularity curves, and the problem of sparsely representing such
singularities using 2D representation systems has been extensively
studied; only curvelets \cite{CD04}, contourlets \cite{DV05}, and
shearlets \cite{GL07} are known to succeed in this task in an optimal
way (see also Section~\ref{sec:sparse-appr}). We describe
contourlets and curvelets in more details in
Section~\ref{sec:applications} and will here just mention some
differences to shearlets. Contourlets are constructed from a
discrete filter bank and have therefore, unlike shearlets, no
continuum theory. Curvelets, on the other hand, are a
continuum-domain system which, unlike shearlets, does not transfer
in a uniform way to the digital world. It is fair to say that
shearlet theory is a comprehensive theory with a mathematically rich
structure as well as a superior connection between the continuum and digital
realm.

The missing link between the continuum and digital world for curvelets
is caused by the use of rotation as a means to parameterize
directions. One of the distinctive features of shearlets is the
use of shearing in place of rotation; this is, in fact,
decisive for a clear link between the continuum and digital world
which stems from the fact that the shear matrix preserves the integer
lattice. Traditionally, the shear parameter ranges over a non-bounded
interval. This has the effect that the directions are not treated
uniformly, which is particularly important in applications. On the
other hand, rotations clearly do not suffer from this deficiency. To
overcome this shortcoming of shearing, Guo, Labate, and Weiss together with
two of the authors \cite{LLKW05} (see also \cite{GKL06}) \gggk{introduced}
the so-called cone-adapted shearlet systems, where the frequency plane is
partitioned into a horizontal and a vertical cone which allows
restriction of the shear parameter to bounded intervals
(Section~\ref{sec:prel-term}), thereby guaranteeing uniform
treatment of directions.

Shearlet systems therefore come in two ways: One class being generated
by a unitary representation of the shearlet group and equipped with a
particularly `nice' mathematical structure, however causes a
bias towards one direction, which makes it unattractive for
applications; the other class being generated by a quite similar
procedure, but restricted to cones in frequency domain, thereby
ensuring an equal treatment of all directions. To be precise this treatment of directions is
only `almost equal' since there still is a slight,
but controllable, bias towards directions of the coordinate
axes, see also Figure~\ref{fig:2d-tiling} in
Section~\ref{sec:class-constr}. For both classes, the
\emph{continuous} shearlet systems are associated with a 4-dimensional
parameter space consisting of a scale parameter measuring the
resolution, a shear parameter measuring the orientation, and a
translation parameter measuring the position of the shearlet (Section
\ref{sec:cont-shearl-transf}). A sampling of this parameter space
leads to \emph{discrete} shearlet systems, and it is obvious that the
possibilities for this are numerous. Using dyadic sampling leads to
so-called regular shearlet systems which are those discrete systems
mainly considered in this paper. It should be mentioned that also
irregular shearlet systems have attracted some attention, and we refer
to the papers \cite{KL07,KLP10,KKL10}. We end this section by
remarking that these discrete shearlet systems belong to a larger
class of representation systems -- the so-called composite wavelets
\cite{GLLWW05, GLLWW06, GLLWW06a}.

\subsection{Anisotropic Features, Discrete Shearlet Systems, and Quest for Sparse Approximations}
\label{sec:anis-feat-sparse}

In many applications in 2D and 3D imaging the important information
is often located around \emph{edges} \gggk{separating} `image
objects' \gggk{from} `background'. These features correspond precisely to
the anisotropic structures in the data. Two-dimensional shearlet
systems are carefully designed to efficiently encode such
anisotropic features. In order to do this effectively, shearlets are
scaled according to a parabolic scaling law, thereby exhibiting a
spatial footprint of size $2^{-j}$ times $2^{-j/2}$, where $2^j$ is
the (discrete) scale parameter; this should be compared to the size
of wavelet footprints: $2^{-j}$ times $2^{-j}$. These elongated,
scaled needle-like shearlets then parametrize directions by slope
encoded in a shear matrix. As mentioned in the previous section,
such carefully designed shearlets do, in fact, perform optimally
when representing and analyzing anisotropic features in 2D data
(Section \ref{sec:sparse-appr}).

In 3D the situation changes somewhat. While in 2D we `only' have to
handle one type of anisotropic structures, namely curves, in 3D a much
more complex situation can occur, since we find two geometrically very
different anisotropic structures: Curves as one-dimensional features
and surfaces as two-dimensional anisotropic features. Our 3D shearlet
elements in spatial domain will be of size $2^{-j}$ times $2^{-j/2}$
times $2^{-j/2}$ which corresponds to `plate-like' elements as $j \to
\infty$. This indicates that these 3D shearlet systems have been
designed to efficiently capture two-dimensional anisotropic
structures, but neglecting one-dimensional structures. Nonetheless,
surprisingly, these 3D shearlet systems still perform optimally when
representing and analyzing 3D data that contain both curve and surface
singularities (Section~\ref{sec:shearl-high-dimens}).

Of course, before we can talk of optimally sparse approximations, we need
to actually have these 2D and 3D shearlet systems at hand.
Several constructions of discrete band-limited 2D shearlet frames are
already known, see \cite{GKL06,KL07,DKST09,KKL10}.
But since spatial localization of the
analyzing elements of the encoding system is immensely important both
for a precise detection of geometric features as well as for a fast
decomposition algorithm, we will mainly follow the sufficient
conditions for and construction of compactly supported cone-adapted
2D shearlet systems by Kittipoom and two of the authors \cite{KLP10}
(Section~\ref{sec:constr-comp-supp}). These results provide a large
class of separable, compactly supported shearlet systems with good frame
bounds, optimally sparse approximation properties, and associated numerically
stable algorithms.

\subsection{Continuous Shearlet Systems}
\label{sec:cont-shearl-transf}

Discrete shearlet systems are, as mentioned, a sampled version of the
so-called continuous shearlet systems. These continuous shearlets
come, of course, also in two different flavors, and we will briefly
describe these in this section.

\subsubsection{Cone-Adapted Shearlet Systems}
\label{sec:cone-adapt-transf}

Anisotropic features in multivariate data can be modeled in many
different ways. One possibility is the cartoon-like image
class discussed above,
but one can also model such directional singularities through
distributions. One would, for example, model a one-dimensional
anisotropic structure as the delta distribution of a curve. The
so-called \emph{cone-adapted continuous shearlet transform} associated with
\emph{cone-adapted continuous shearlet systems} was
introduced by Labate and the first author in \cite{KL09} in the study
of resolutions of the wavefront set for such distributions. It was shown
that the continuous shearlet transform is not only able to identify
the singular support of a distribution, but also the
\emph{orientation} of distributed singularities along curves.
More precisely, for a class of band-limited shearlet generators $\psi \in
L^2(\R^2)$, the first author and Labate \cite{KL09} showed that the
wavefront set of a (tempered) distribution $f$ is precisely the
closure of the set of points $(t,s)$, where the shearlet transform of
$f$ \[
\gggk{(a,s,t) \mapsto} \innerprod{f}{a^{-3/4}\psi(A_a^{-1}S_s^{-1}(\cdot-t))}, \quad
\text{where } A_a = \begin{pmatrix} a & 0 \\ 0 & a^{1/2}
\end{pmatrix}
\text{ and }
S_s = \begin{pmatrix}
1 & s \\ 0 & 1
\end{pmatrix},
\]
is \emph{not} of fast decay as the scale parameter $a\to 0$. Later
Grohs \cite{Gro09b} extended this result to Schwartz-class generators
with infinitely many directional vanishing moments, in particular, not
necessarily band-limited generators. In other words, these results
demonstrate that the wavefront set of a distribution can be
\emph{precisely captured} by continuous shearlets. For
  constructions of continuous shearlet frames with compact support, we
  refer to \cite{Gro10}.

\subsubsection{Shearlets from Group Representations}
\label{sec:group-transf}

Cone-adapted continuous shearlet systems and their associated
cone-adapted continuous transforms described in the previous section
have only very recently -- in 2009 -- attracted attention.
Historically, the continuous shearlet transform was first introduced
in \cite{GKL06} \gggk{without restriction to cones in frequency domain}.
Later, it was shown in \cite{DKMSST08} that the
associated continuous shearlet systems are generated by a strongly
continuous, irreducible, square-integrable representation of a locally
compact group, the so-called \emph{shearlet group}. This implies that these
shearlet systems possess a rich mathematical structure, which in
\cite{DKMSST08} was \gggk{used to derive} uncertainty principles to
tune the accuracy of the shearlet transform, and which in \cite{DKST09}
\gggk{allowed the usage of} coorbit theory to study smoothness spaces associated
with the decay of the shearlet coefficients.

Dahlke, Steidl, and Teschke generalized the shearlet group and the associated
continuous shearlet transform to higher dimensions $\R^n$ in the paper
\cite{DST09}. Furthermore, in \cite{DST09} they showed that, for
certain band-limited generators, the continuous shearlet transform is
able to identify hyperplane and tetrahedron singularities. Since this transform
originates from a unitary group representation, it is
not able to capture all directions, in particular, it will not capture
the delta distribution on the $x_1$-axis (and more generally, any
singularity with `$x_1$-directions'). We also remark that the
extension in \cite{DST09} uses another scaling matrix as compared to
the one used for the three-dimensional shearlets considered in this
paper; we refer to Section~\ref{sec:shearl-high-dimens} for a more
detailed description of this issue.

\subsection{Applications}
\label{sec:applications}
Shearlet theory has applications in various areas. In this
section we will present two examples of such: Denoising of images and
geometric separation of data. Before, in order to show the reader the
advantages of digital shearlets, we first give a short overview of
the numerical aspects of shearlets and two similar implementations of
directional representation systems, namely contourlets and curvelets,
discussed in Section~\ref{sec:direct-repr-syst}.

\begin{figure}[ht!]
\vspace{-1em}
\centering
\subfloat[Original image]{%
\includegraphics[scale=0.30]{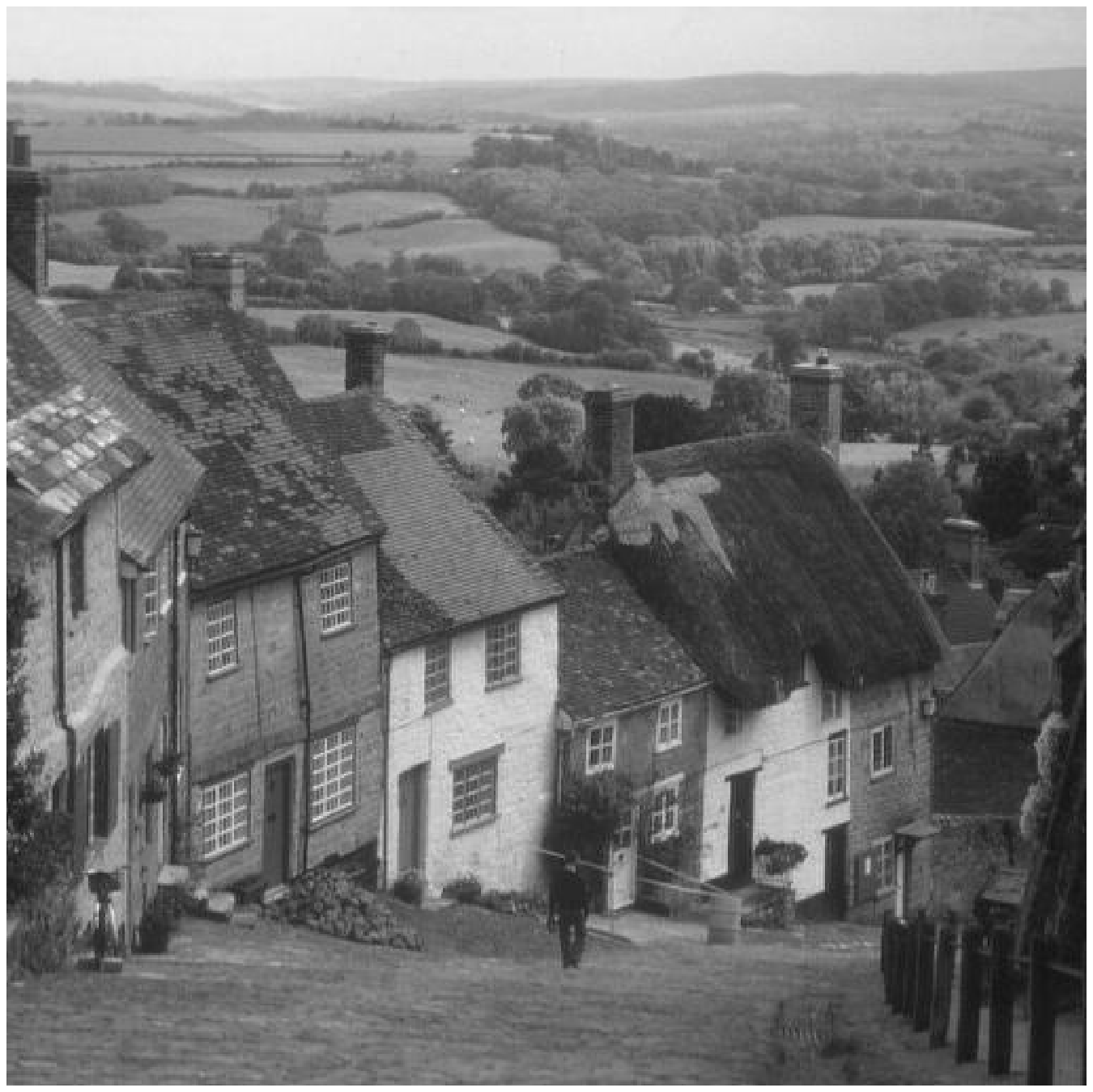}
\label{fig:goldhill_orig}}
\subfloat[Noisy image]{%
\includegraphics[scale=0.30]{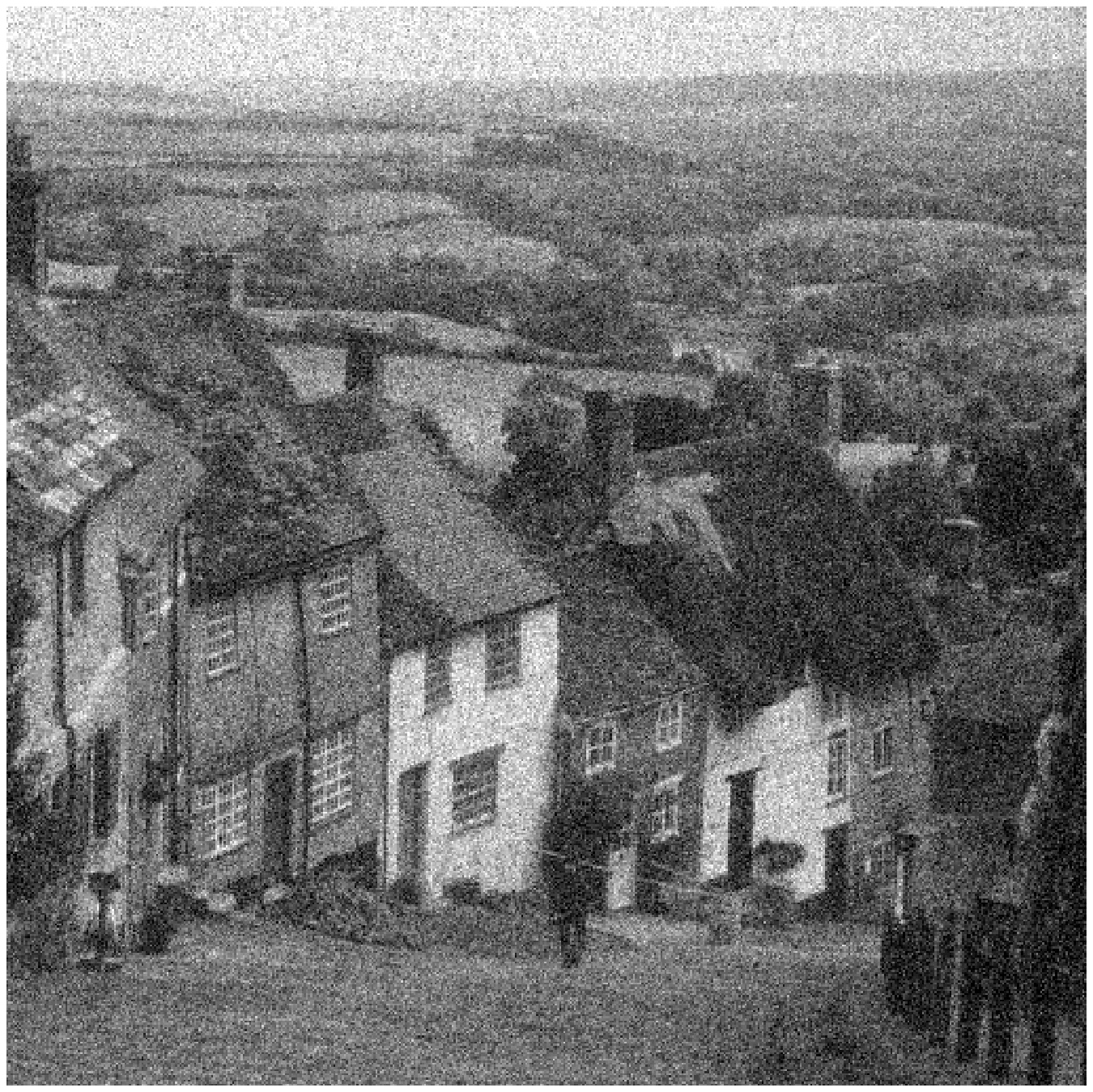}
\label{fig:noisy}}

\mbox{\subfloat[Denoised using curvelets]{%
\includegraphics[scale=0.30]{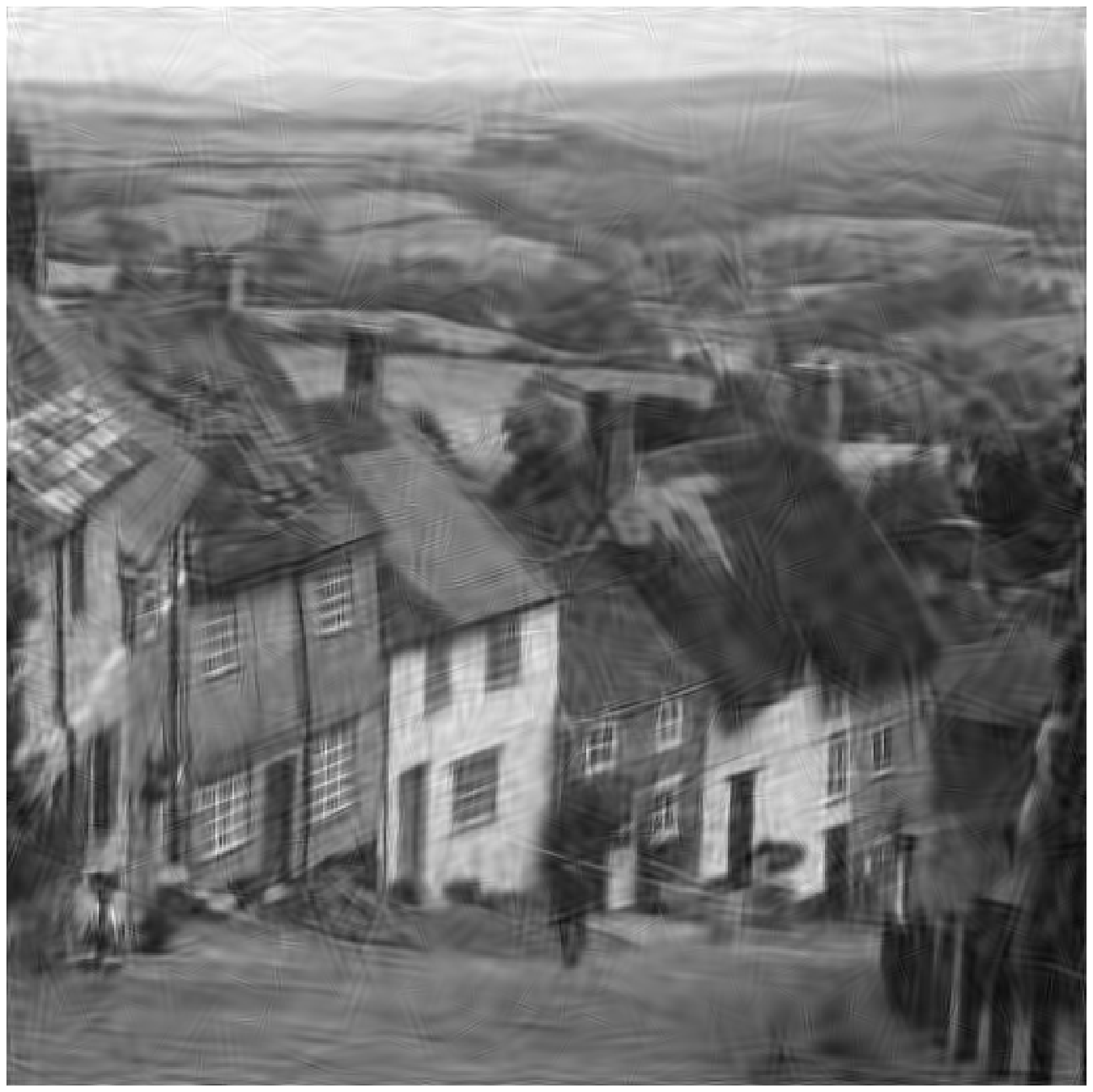}
\label{fig:deno_curve}
}
\subfloat[Denoised using shearlets]{%
\includegraphics[scale=0.30]{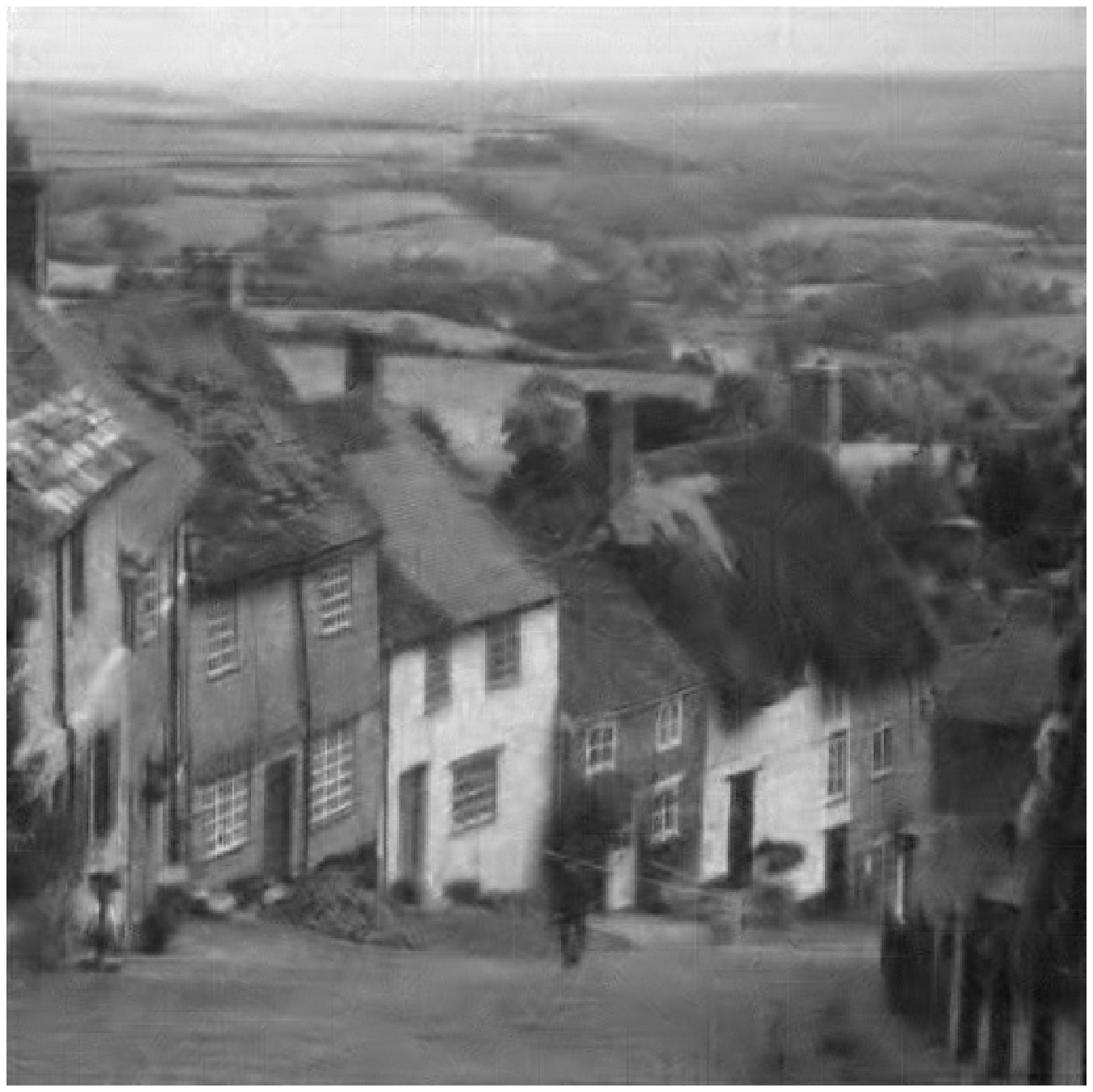}
\label{fig:deno_shear}
}}

\subfloat[Denoised using curvelets (zoom)]{%
\includegraphics[scale=0.30]{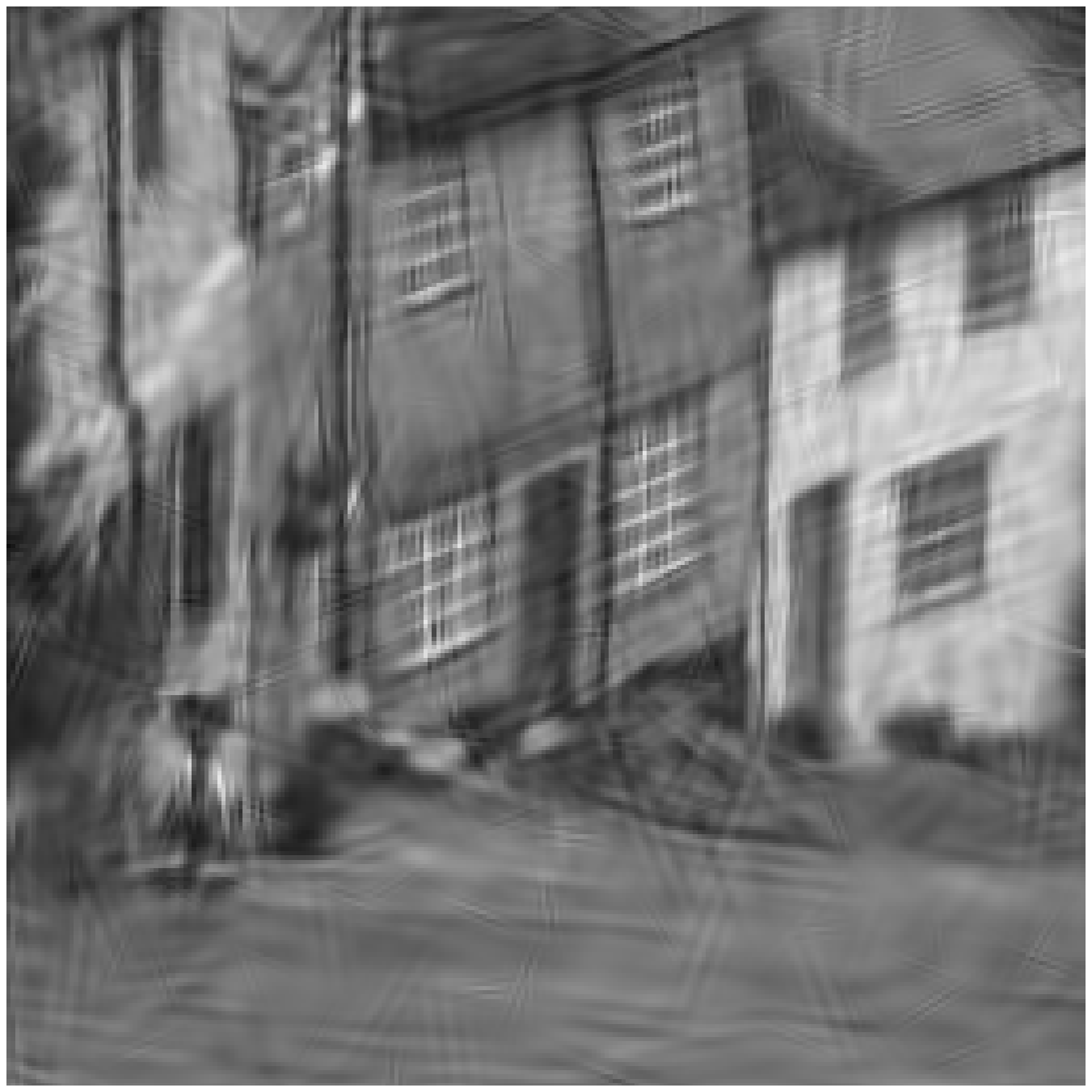}
\label{fig:deno_curve_zoom}
}
\subfloat[Denoised using shearlets (zoom)]{%
\includegraphics[scale=0.30]{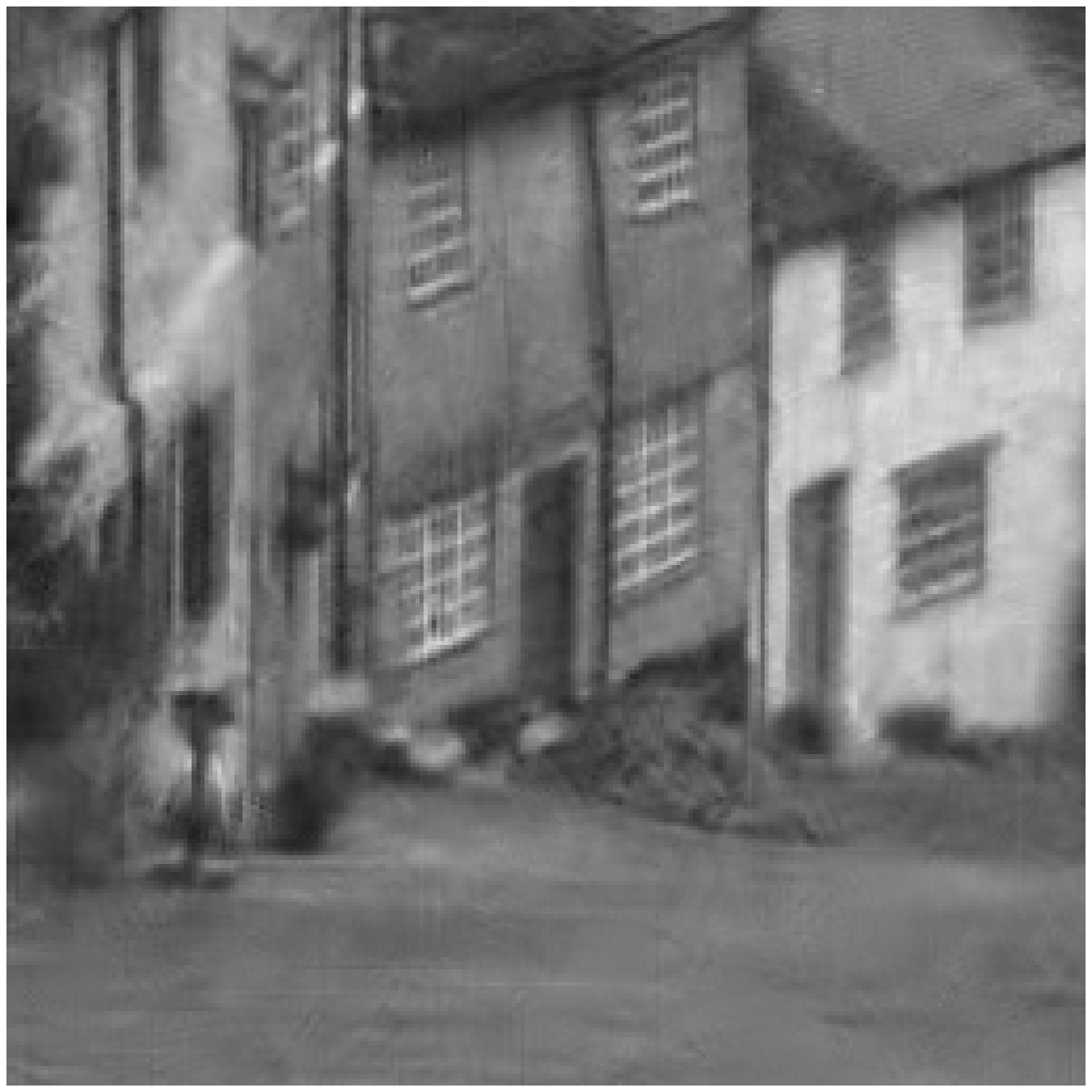}
\label{fig:deno_shear_zoom}
}
\caption{Denoising of the Goldhill image ($512 \times 512$) using
  shearlets and curvelets. The noisy image in (b) has a peak
    signal-to-noise ratio of $20.17$ dB. The curvelet-denoised image
    in (c) and~(e) has a PSNR of $28.70$ dB, while the shearlet-denoised
    image in (d) and~(f) has a PSNR of only $29.20$ dB. }
\label{fig:denoising}
\end{figure}
\begin{description}
\item[{\em Curvelets} \cite{CDDY06}.] \hspace*{-0.4cm} This approach builds on directional
  frequency partitioning and the use of the Fast Fourier transform.
  The algorithm can be efficiently implemented using (in frequency
  domain) multiplication with the frequency response of a filter and
  frequency wrapping in place of convolution and down-sampling.
  However, curvelets need to be band-limited and can only have very
  good spatial localization if one allows high redundancy.
\item[{\em Contourlets} \cite{DV05}.] \hspace*{-0.4cm}
  This approach uses a directional filter bank, which produces
  directional frequency partitioning similar to those of curvelets. As
  the main advantage of this approach, it allows a tree-structured
  filter bank implementation, in which aliasing due to subsampling is
  allowed to exist. Consequently, one can achieve great  efficiency  in
  terms of redundancy and good spatial localization. However, the
  directional selectivity in this approach is artificially imposed by
  the special sampling rule of a filter bank which introduces various
  artifacts. We remark that also the recently introduced {\em Hybrid Wavelets} \cite{ER07}
  suffer from this deficiency.
\item[{\em Shearlets} \cite{Lim09b}.] \hspace*{-0.4cm}  Using a shear matrix instead of
  rotation, directionality is naturally adapted for the digital
  setting in the sense that the shear matrix preserves the structure
  of the integer grid. Furthermore, excellent spatial localization is
  achieved by using compactly supported shearlets. The only
  drawback is that these compactly supported shearlets are not tight
  frames and, accordingly, the synthesis process needs to be performed
  by iterative methods.
\end{description}

To illustrate how two of these implementations perform, we have
included a denoising example of the Goldhill image using both
curvelets\footnote{Produced using Curvelab (Version 2.1.2), which is
  available from \url{http://curvelet.org}.} and shearlets, see
Figure~\ref{fig:denoising}. We omit a detailed analysis of the
denoising results and leave the visual comparison to the reader. For
a detailed review of the shearlet transform and associated aspects,
we refer to \cite{Lim09b,DKS08,ELL08, DKSZ10}. We also refer to
\cite{KS10,hks10} for MRA based algorithmic approaches to the shearlet
transform.

\begin{figure}
\centering
\subfloat[Original point-curve data of size $256 \times 256$.]{%
\includegraphics[width=3.5cm,height=3.5cm]{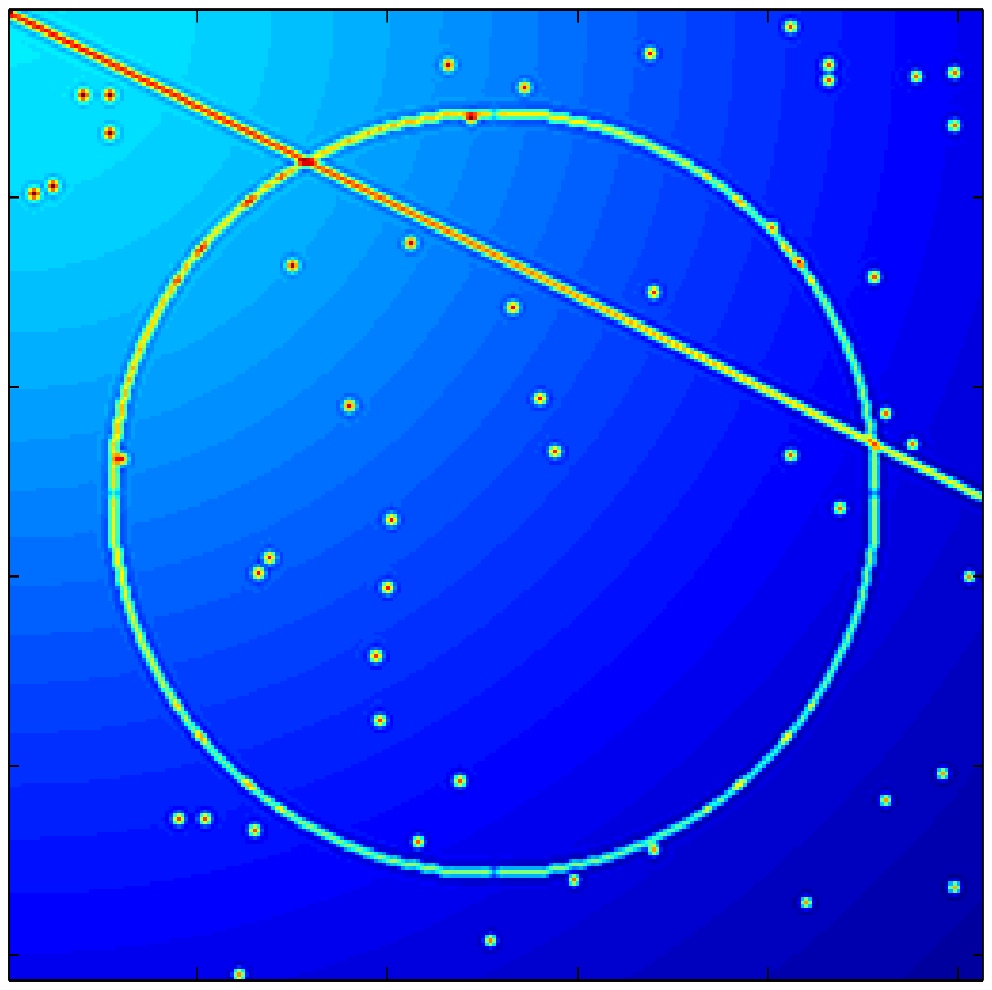}
\label{fig:pointcurve}
}
\quad
\subfloat[Separated point-like data (captured by wavelets).]{%
\includegraphics[width=3.5cm,height=3.5cm]{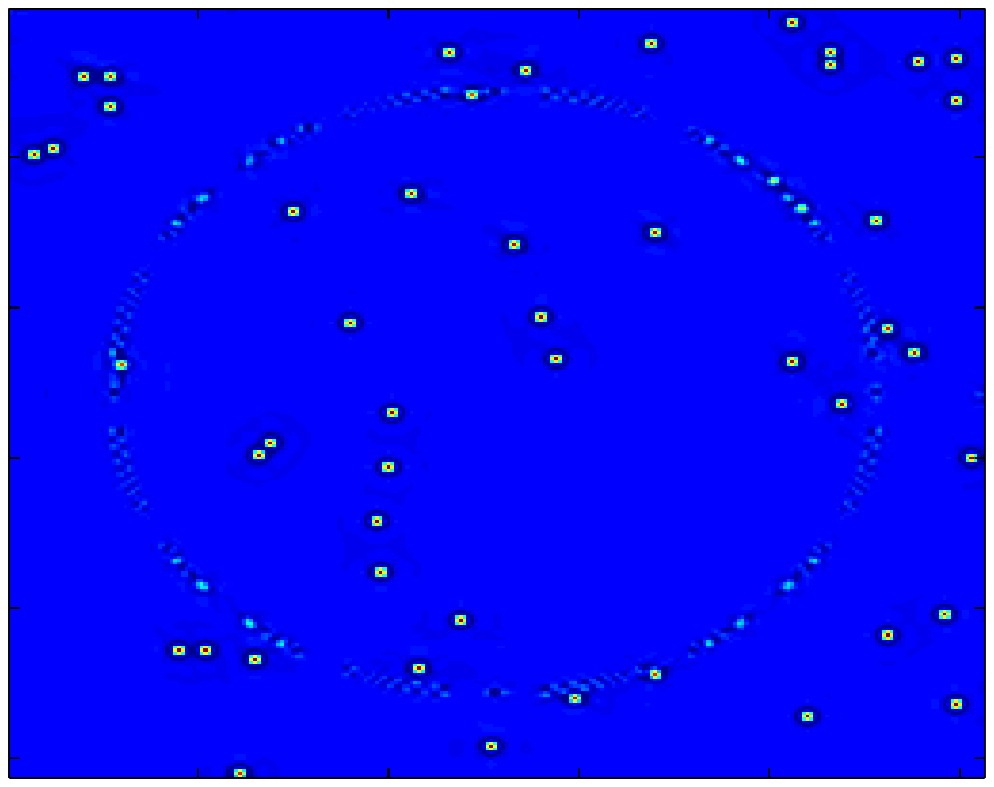}
\label{fig:points}
}
\quad
\subfloat[Separated curve-like data (captured by shearlets).]{%
\includegraphics[width=3.5cm,height=3.5cm]{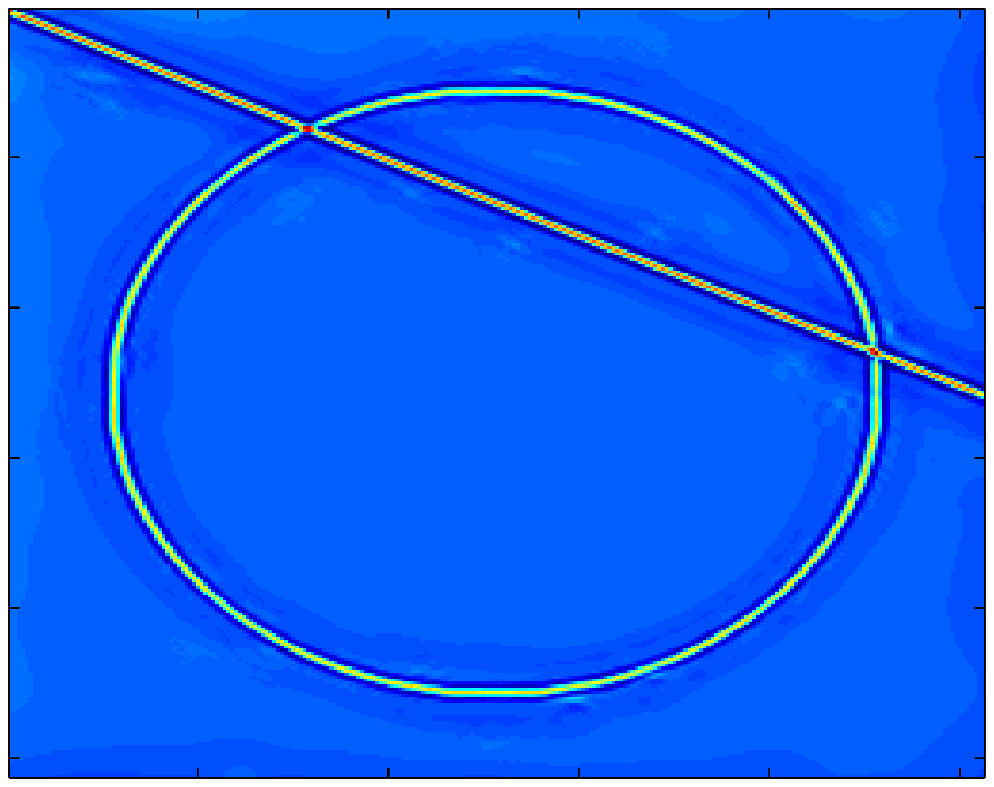}
\label{fig:curves}
}
\caption{Geometric separation of mixed `point-and-curve' data.
(a): Input data. (b) and (c): The output of the separation algorithm.}
\label{fig:geosep}
\end{figure}

The shearlet transform, in companion with the wavelet transform, has
also been applied to accomplish geometric separation of
`point-and-curve'-like data. An artificially made example of such
data can be seen in Figure~\ref{fig:pointcurve}. For a theoretical
account of these separation ideas we refer to the recent papers by
Donoho and the first author \cite{DK09,DK10}. Here we simply display
the result of the separation, see Figure~\ref{fig:geosep}. For
real-world applications of these separation techniques we refer to
the paper \cite{avignon} on neurobiological imaging.

In the spirit of reproducible research \cite{DMSSU09}, we wish
to mention that Figure~\ref{fig:deno_shear}, \ref{fig:deno_shear_zoom}
and \ref{fig:geosep} have been produced by the discrete shearlet
transform implemented
in the Matlab toolbox \emph{Shearlab}
which has recently been released under
a GNU license and is freely available at
\url{http://www.shearlab.org}.

\subsection{Outline}
\label{sec:outline}
In Section~\ref{sec:shearl-two-dimens} we present a review of shearlet
theory in $L^2(\R^2)$, where we focus on discrete shearlet systems. We describe the classical band-limited
construction (Section~\ref{sec:class-constr}) and a more recent
construction of compactly supported shearlets
(Section~\ref{sec:constr-comp-supp}). In Section~\ref{sec:sparse-appr}
we present results on the ability of shearlets to optimally sparsely
approximate cartoon-like images.
Section~\ref{sec:shearl-high-dimens} is dedicated to a discussion
on similar properties of 3D shearlet systems.

\section{2D Shearlets}
\label{sec:shearl-two-dimens}
In this section, we summarize what is known about constructions of
discrete shearlet systems in 2D. Although all results in this section can easily be
extended to (irregular) shearlet systems associated with a general
irregular set of parameters for scaling, shear, and translation, we
will only focus on the discrete shearlet systems associated with a
regular set of parameters as described in the next section. For a
detailed analysis of irregular shearlet systems, we refer to
\cite{KLP10}. We first start with various notations and definitions
for later use.
\subsection{Preliminaries}
\label{sec:prel-term}
For $j \ge 0, k \in \mathbb{Z}$, let
\[
A_{2^j} = \begin{pmatrix}
2^j & 0 \\ 0 & 2^{j/2}
\end{pmatrix}
,\,\,
S_k = \begin{pmatrix}
1 & k \\ 0 & 1
\end{pmatrix},
\quad \text{and} \quad M_c =
\begin{pmatrix}
c_1 & 0 \\ 0 & c_2
\end{pmatrix},
\]
where $c = (c_1,c_2)$  and $c_1,c_2$ are some positive constants.
Similarly, we define
\[
\tilde{A}_{2^j} = \begin{pmatrix}
2^{j/2} & 0 \\ 0 & 2^{j}
\end{pmatrix}
,\,\,
\tilde{S}_k = \begin{pmatrix}
1 & 0 \\ k & 1
\end{pmatrix},
\quad \text{and} \quad \tilde{M}_c =
\begin{pmatrix}
c_2 & 0 \\ 0 & c_1
\end{pmatrix}.
\]
Next we define discrete shearlet systems in 2D.
\begin{definition}
  Let $c=(c_1,c_2) \in (\R_+)^2$. For $\phi, \psi, \tilde{\psi} \in
  L^2(\R^2)$ the \emph{cone-adapted 2D discrete shearlet system}
  $\SH(\phi,\psi,\tilde{\psi};c)$ is defined by
\[
\SH(\phi,\psi,\tilde{\psi};c) = \Phi(\phi;c_1) \cup \Psi(\psi;c) \cup
\tilde{\Psi}(\tilde{\psi};c),
\]
where
\begin{align*}
  \Phi(\phi;c_1) &= \{\phi(\cdot-m) : m \in c_1\Z^2\},\\
  \Psi(\psi;c) &= \{2^{\frac34 j} \psi(S_kA_{2^j}\, \cdot \, -m) : j \ge 0,
  -\lceil2^{j/2}\rceil \leq k \leq \lceil2^{j/2}\rceil, m \in M_c\Z^2
  \},
  \intertext{and}
  \tilde{\Psi}(\tilde{\psi};c) &= \{2^{\frac34 j}
  \tilde{\psi}(\tilde{S}_k\tilde{A}_{2^j}\, \cdot \, -m) : j \ge 0,
  -\lceil 2^{j/2} \rceil \leq k \leq \lceil 2^{j/2} \rceil, m \in
  \tilde{M}_c\Z^2 \}.
\end{align*}
If $\SH(\phi,\psi,\tilde{\psi};c)$ is a frame for $L^2(\R^2)$, we
refer to $\phi$ as a \emph{scaling function} and $\psi$ and
$\tilde{\psi}$ as \emph{shearlets}.
\end{definition}
Our aim is to construct compactly supported functions $\phi,\psi$, and
$\tilde{\psi}$ to obtain compactly supported shearlets in 2D. For
this, we will describe general sufficient conditions \jl{on the}
shearlet generators $\psi$ and $\tilde{\psi}$, which lead to the
construction of compactly supported shearlets. To formulate our
sufficient conditions on $\psi$ and $\tilde{\psi}$
(Section~\ref{sec:constr-comp-supp}), we will first need to introduce the
necessary notational concepts.

For functions $\phi,\psi,\tilde{\psi} \in L^2(\R^2)$, we define
$\Theta : \R^2 \times \R^2 \to \R$  by
\begin{equation}\label{eq:moPhi}
\Theta(\xi,\omega) = |\hat \phi(\xi)||\hat \phi(\xi+\omega)|+ \Theta_1(\xi,\omega)+\Theta_2(\xi,\omega),
\end{equation}
where
\[
\Theta_1(\xi,\omega) = \sum\limits_{j \ge 0} \sum\limits_{|k| \le \lceil 2^{j/2} \rceil}
\left|\hat{\psi}(S_{{k}}^{T} A_{{2^{-j}}}\xi) \right| \left|\hat{\psi}({S_{k}}^{T}A_{{2^{-j}}}\xi + \omega) \right|
\]
and
\[
\Theta_2(\xi,\omega)=\sum\limits_{j \ge 0} \sum\limits_{|k| \le \lceil 2^{j/2} \rceil}
\left|\hat{\tilde{\psi}}({S}_{{k}} \tilde{A}_{{2^{-j}}}\xi) \right|
\left|\hat{\tilde{\psi}}({S}_{{k}}\tilde{A}_{{2^{-j}}}\xi+ \omega) \right|.
\]
Also, for $c = (c_1,c_2) \in (\R_+)^2$, let
\begin{align*}
R(c) &= \sum_{m \in \Z^2 \setminus \{0\}} \left(\Gamma_0(c_1^{-1}m)\Gamma_0(-c_1^{-1}m)\right)^{\frac 12}
+\left(\Gamma_1(M_c^{-1}m)\Gamma_1(-M_c^{-1}m)\right)^{\frac 12} \\
&\phantom{\sum_{m \in \Z^2 \setminus \{0\}} +} + (\Gamma_2(\tilde{M}_c^{-1}m)\Gamma_2(-\tilde{M}_c^{-1}m))^{\frac 12},
\end{align*}
where
\[
\Gamma_0(\omega) = \esssup_{\xi \in \R^2} |\hat \phi(\xi)||\hat \phi(\xi+\omega)| \quad \text{and} \quad
\Gamma_{i}(\omega) = \esssup_{\xi \in \R^2} \Theta_{i}(\xi,\omega) \quad \text{for}\,\, i = 1,2.
\]

\subsection{Classical Construction}
\label{sec:class-constr}
We now first describe the construction of band-limited
shearlets which provides tight frames for $L^2(\R^2)$. Constructions of
this type were first introduced by Labate, Weiss, and two of the authors
in \cite{LLKW05}. The \emph{classical example} of a generating
shearlet is a function $\psi \in L^2(\R^2)$ satisfying
\[
\hat{\psi}(\xi) = \hat{\psi}(\xi_1,\xi_2) = \hat{\psi}_1(\xi_1) \, \hat{\psi}_2(\tfrac{\xi_2}{\xi_1}),
\]
where $\psi_1 \in L^2(\R)$ is a discrete wavelet, \ie satisfies the
discrete Calder\'{o}n condition given by
\[
\sum_{j \in \Z}\abssmall{\hat\psi_1(2^{-j}\xi)}^2 = 1 \quad \mbox{for a.e. } \xi \in \R,
\]
with $\hat{\psi}_1 \in C^\infty(\R)$ and  $\supp \hat{\psi}_1
\subseteq [-\frac54,-\frac14] \cup [\frac14,\frac54]$, and $\psi_2 \in
L^2(\R)$ is a bump function, namely
\[
\sum_{k = -1}^{1}\abssmall{\hat\psi_2(\xi+k)}^2 = 1 \quad \mbox{for a.e. } \xi \in [-1,1],
\]
satisfying $\hat{\psi}_2 \in C^\infty(\R)$ and $\supp \hat{\psi}_2 \subseteq \itvcc{-1}{1}$. There are several choices of
$\psi_{1}$ and $\psi_{2}$ satisfying those conditions, and we refer to
\cite{GKL06} for further details.
\begin{figure}%
\centering
\parbox[t]{0.35\textwidth}{%
\centering
\includegraphics[height=1.4in]{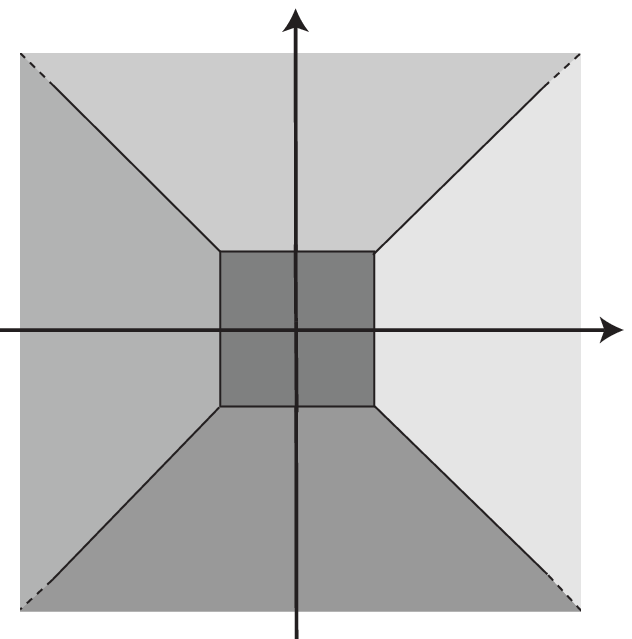}
\put(-33,58){\footnotesize{$\cC_1$}}
\put(-70,80){\footnotesize{$\cC_2$}}
\put(-88,30){\footnotesize{$\cC_3$}}
\put(-50,52){\footnotesize{$\cR$}}
\put(-45,15){\footnotesize{$\cC_4$}}
\caption{The cones $\cC_1$ -- $\cC_4$ and the centered rectangle $\cR$
  in the frequency domain.}%
\label{fig:cones-rectangle}}%
 \hspace*{.15\textwidth}
\begin{minipage}[t]{0.33\textwidth}%
\centering
\includegraphics[height=1.4in]{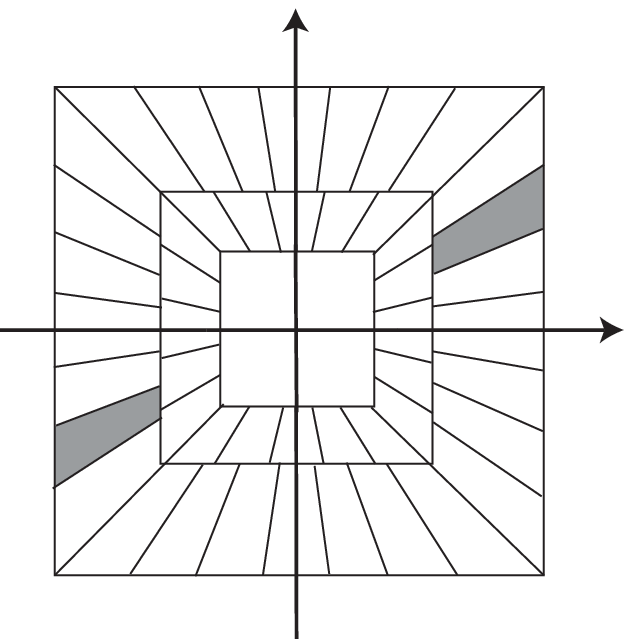}
\caption{Tiling of the frequency domain induced by band-limited shearlets.}%
\label{fig:2d-tiling}%
\end{minipage}%
\end{figure}%
The tiling of the frequency domain given by these band-limited
generators and choosing $\tilde{\psi}(x_1,x_2)=\psi(x_2,x_1)$ is
illustrated in Figure~\ref{fig:2d-tiling}. As described in
Figure~\ref{fig:cones-rectangle}, a conic region $\cC_1 \cup \cC_3$ is
covered by the frequency support of shearlets in $\Psi(\psi;c)$ while
$\cC_2 \cup \cC_4$ is covered by $\tilde\Psi(\tilde \psi;c)$. For this
particular choice, using an appropriate scaling function $\phi$ for
the centered rectangle $\cR$ (see Figure~\ref{fig:cones-rectangle}),
it was proved in \cite[Thm.\ 3]{GKL06} that the associated
cone-adapted discrete shearlet system
$\SH(\phi,\psi,\tilde{\psi};(1,1))$ forms a Parseval frame for
$L^2(\R^2)$.

\subsection{Constructing Compactly Supported Shearlets}
\label{sec:constr-comp-supp}
We are now ready to state general sufficient conditions for the
construction of shearlet frames.

\begin{theorem}[\protect\hspace*{-.248em}\cite{KLP10}]
\label{theo:general}
Let $\phi, \psi \in L^2(\R^2)$ be functions such that
\[\hat \phi(\xi_1,\xi_2) \leq C_1 \cdot \min{\{1,\abs{\xi_1}^{-\gamma}\}} \cdot \min{\{1,|\xi_2|^{-\gamma}\}}\]
and
\begin{equation}\label{eq:decay}
\abs{\hat{\psi}(\xi_{1}, \xi_{2})} \le C_2 \cdot \min\{1,|\xi_1|^{\alpha}\} \cdot \min{\{1,|\xi_1|^{-\gamma}\}} \cdot \min{\{1,|\xi_2|^{-\gamma}\}},
\end{equation}
for some positive constants $C_1,C_2 < \infty$ and $\alpha>\gamma>3$.
Define $\tilde \psi(x_1,x_2) = \psi(x_2,x_1)$, and let $L_{\inf},
L_{\sup}$ be defined by
\[
L_{\inf} = \essinf\limits_{\xi \in \R^2} \wql{\Theta(\xi,0)} \quad \text{and} \quad L_{\sup} = \esssup\limits_{\xi \in \R^2} \wql{\Theta(\xi,0)}.
\]
Suppose that there is a constant $\tilde{L}_{\inf}>0$ such that
$0<\tilde{L}_{\inf}\le L_{\inf}$. Then there exist a sampling
parameter $c = (c_1,c_2)$ with $c_1 = c_2$ and a constant
$\tilde{L}_{sup} < \infty$ such that
\[
R(c) < \tilde{L}_{inf} \le L_{inf} \,\, \mbox{and} \,\, L_{sup} \le \tilde{L}_{sup},
\]
and, further, ${\SH}(\phi,\psi,\tilde \psi;c)$ forms a frame for
$L^2(\R^2)$ with frame bounds $A$ and $B$ satisfying
\begin{equation}\label{eq:main_2d}
\frac{1}{|\det M_c|} [\tilde{L}_{inf} - R(c)] \le A \le B \le \frac{1}{|\det M_c|} [\tilde{L}_{sup} + R(c)].
\end{equation}
\end{theorem}
For a detailed proof, we refer to the paper \cite{KLP10} by Kittipoom
and two of the authors.

Obviously, band-limited shearlets (from
Section~\ref{sec:class-constr}) satisfy condition \eqref{eq:decay}.
More interestingly, also a large class of spatially compactly
supported functions satisfies this condition. In fact, in
\cite{KLP10}, various constructions of compactly supported shearlets
are presented using Theorem \ref{theo:general} and generalized
low-pass filters; an example of such a construction procedure is given
in Theorem~\ref{thm:low-pass-construction} below. In
Theorem~\ref{theo:general} we assumed $c_1 = c_2$ for the sampling
matrix $M_c$ (or $\tilde{M}_c$), the only reason for this being the
simplification of the estimates for the frame bounds $A,B$ in
\eqref{eq:main_2d}. In fact, the estimate
\eqref{eq:main_2d} 
generalizes easily to non-uniform sampling constants $c_1,c_2$ with
$c_1 \neq c_2$. For explicit estimates of the form \eqref{eq:main_2d}
in the case of non-uniform sampling, we refer to \cite{KLP10}.

The following result provides a specific family of functions satisfying
the general sufficiency condition from Theorem~\ref{theo:general}.
\begin{theorem}[\protect\hspace*{-.255em}\ \cite{KLP10}]
\label{thm:low-pass-construction}
Let $K, L > 0$ be such that $L \ge 10$ and $\frac{3L}{2} \le K \le 3L-2$, and define a shearlet $\psi \in L^2(\R^2)$ by
\[
\hat{\psi}(\xi) = m_1(4\xi_1)\hat{\phi}(\xi_1)\hat{\phi}(2\xi_2), \quad \xi = (\xi_1,\xi_2) \in \R^2,
\]
where $m_0$ is the low pass filter satisfying
\[
|m_0(\xi_1)|^2 = (\cos(\pi\xi_1))^{2K}\sum_{n=0}^{L-1} \binom{K-1+n}{ n}(\sin(\pi\xi_1))^{2n}, \quad \xi_1 \in \R,
\]
$m_1$ is the associated bandpass filter defined by
\[
|m_1(\xi_1)|^2 = |m_0(\xi_1+\tfrac12)|^2, \quad \xi_1 \in \R,
\]
and $\phi$ is the scaling function given by
\[
\hat{\phi}(\xi_1) = \prod_{j=0}^{\infty} m_0(2^{-j}\xi_1), \quad \xi_1 \in \R.
\]
Then there exists a sampling constant $\hat c_1>0$ such that the
shearlet system $\Psi(\psi;c)$ forms a frame for $\check L^2(\cC_1
\cup \cC_3):=\setprop{f \in L^2(\R^2)}{\supp \hat f \subset \cC_1 \cup
  \cC_3}$ for any sampling matrix $M_c$ with $c=(c_1,c_2) \in
(\R_+)^2$ and $c_2 \le c_1 \le \hat{c}_1$.
\end{theorem}

For these shearlet systems,
there is a bias towards the vertical axis, especially at coarse scales, since
they are defined for $\check L^2(\cC_1 \cup \cC_3)$, and hence, the
frequency support of the shearlet elements overlaps more significantly
along the vertical axis.
In order to control the
upper frame bound, it is therefore desirable to apply a denser
sampling along the vertical axis than along the horizontal axis, i.e.,
$c_1 > c_2$.

Having compactly supported (separable) shearlet frames for $\check
L^2(\cC_1 \cup \cC_3)$ at hand by
Theorem~\ref{thm:low-pass-construction}, we can easily construct
shearlet frames for the whole space $L^2(\R^2)$. The exact procedure
is described in the following theorem from \cite{KLP10}.

\begin{theorem}[\hspace*{-.23em}\cite{KLP10}]
 \label{thm:frame-for-R2}
 Let $\psi \in L^2(\R^2)$ be the shearlet with associated scaling
 function $\gggk{\phi_1} \in L^2(\R)$ both introduced in
 Theorem~\ref{thm:low-pass-construction}, and set
 $\phi(x_1,x_2)=\gggk{\phi_1(x_1) \phi_1(x_2)}$ and
 $\tilde{\psi}(x_1,x_2) = \psi(x_2,x_1)$. Then the corresponding
 shearlet system $\SH(\phi,\psi,\tilde{\psi};c)$ forms a frame for
 $L^2(\R^2)$ for any sampling matrices $M_c$ and $\tilde{M}_c$ with
 $c=(c_1,c_2) \in (\R_+)^2$ and $c_2 \le c_1 \le \hat{c}_1$.
\end{theorem}
For the horizontal cone $\cC_1 \cup \cC_3$ we allow for a denser
sampling by $M_c$ along the vertical axis, i.e., $c_2 \le c_1$,
precisely as in Theorem~\ref{thm:low-pass-construction}. For the
vertical cone $\cC_2 \cup \cC_4$ we analogously allow for a denser
sampling along the horizontal axis; since the position of $c_1$ and
$c_2$ is reversed in $\tilde{M}_c$ compared to $M_c$, this still
corresponds to $c_2 \le c_1$.

We wish to mention that there is a trade-off between \emph{compact
  support} of the shearlet generators, \emph{tightness} of the
associated frame, and \emph{separability} of the shearlet generators.
The known constructions of tight shearlet frames do not use separable
generators (Section~\ref{sec:class-constr}), and these constructions
can be shown to \emph{not} be applicable to compactly supported
generators. Tightness is difficult to obtain while allowing for
compactly supported generators, but we can gain separability as in
Theorem~\ref{thm:frame-for-R2}, hence fast algorithmic realizations.
On the other hand, when allowing non-compactly supported generators,
tightness is possible, but separability seems to be out of reach,
which makes fast algorithmic realizations very difficult.

We end this section by remarking that the construction results above
even generalize to constructions of irregular shearlet systems
\cite{KL07,KKL10}.

\section{Sparse Approximations}
\label{sec:sparse-appr}
After having introduced compactly supported shearlet systems in the previous
section, we  now aim for optimally sparse approximations. To be
precise, we will show that these compactly supported shearlet systems
provide optimally sparse approximations  when representing and
analyzing anisotropic features in 2D data.

\subsection{Cartoon-Like Image Model}
\label{sec:cartoon-like-image}
Following \cite{Don01}, we introduce $STAR^2(\nu)$, a class of
sets $B$ with $C^2$ boundaries $\partial B$ and
curvature bounded by $\nu$, as well as $\cE^2_\nu(\R^2)$, a class of
cartoon-like images. For this, in polar coordinates, we let
$\rho\jll{: \itvco{0}{2\pi}} \rightarrow [0,1]$ be a radius function and define the
set $B$ by
\[
B = \{x \in \R^2 : \abs{x} \le \rho(\theta), \, x = (\abs{x},\theta) \text{ in polar coordinates}\}.
\]
In particular, \jll{we will require that} the boundary $\partial B$ of $B$ is given by the curve
\begin{equation}\label{eq:curve}
\beta(\theta) = \begin{pmatrix} \rho(\theta)\cos(\theta) \\
\rho(\theta)\sin(\theta)\end{pmatrix},
\end{equation}
 and the class of boundaries of interest to us are defined by
\begin{equation}\label{eq:curvebound}
\sup\abssmall{\rho^{''}(\theta)} \leq \nu, \quad \rho \leq \rho_0 < 1,
\end{equation}
where $\rho_0<1$ needs to be chosen so that $y+B \subset \itvcc{0}{1}^2$ for
some $y\in \R^2$.

The following definition now introduces a class of cartoon-like images.

\begin{definition}
\label{def:cartoon-like-image}
For $\nu > 0$, the set $STAR^2(\nu)$ is defined to be the set of all
$B \subset [0,1]^2$ such that $B$ is a translate of a set obeying
\eqref{eq:curve} and \eqref{eq:curvebound}. Further, $\cE^2_\nu(\R^2)$
denotes the set of functions $f \in L^2(\R^2)$ of the form
\[
f = f_0 + f_1 \chi_{B},
\]
where $B \in STAR^2(\nu)$ and $f_0,f_1 \in C_0^2(\R^2)$ with  $\supp f_i
\subset \itvcc{0}{1}^2$ and $\|f_i\|_{C^2} = \sum_{|\alpha| \leq 2}
\|D^{\alpha}f_i\|_{\infty} \leq 1$ for $ i=1,2$.
\end{definition}
One can also consider a more sophisticated class of cartoon-like
images, where the boundary of $B$ is allowed to be \emph{piecewise}
$C^2$, and we refer to the recent paper by two of the authors
\cite{KL10SoBD} and to similar considerations for the 3D case in
Section~\ref{sec:sparse-appr-3d}.

Donoho \cite{Don01} proved that the optimal approximation rate for
such cartoon-like image models $f \in \cE^2_\nu(\R^2)$ which can be
achieved for almost any representation system under a so-called
polynomial depth search selection procedure of the \gggk{selected
  system elements} is
\[
\norm[2]{f-f_N}^2 \le C \cdot N^{-2} \quad \text{as } N \to \infty,
\]
where $f_N$ is the best $N$-term approximation of $f$. As discussed in
the next section shearlets in 2D do indeed deliver this optimal
approximation rate.

\subsection{Optimally Sparse Approximation of Cartoon-Like Images}
\label{sec:optim-sparse-appr}

Let $\SH(\phi,\psi,\tilde{\psi};c)$  be a shearlet frame for
$L^2(\R^2)$. Since this is a countable set of functions, we can
denote it by $\SH(\phi,\psi,\tilde{\psi};c) = (\sigma_i)_{i \in I}$.
  We let  $(\tilde{\sigma}_i)_{i \in I}$ be a dual frame of
$(\sigma_i)_{i \in I}$.  As our $N$-term approximation $f_N$
of a cartoon-like image $f \in \cE^2_\nu(\R^2)$ by the frame
$\SH(\phi,\psi,\tilde{\psi};c)$, we then take
\[
f_N = \sum_{i \in I_N} \innerprod{f}{\sigma_i}\tilde{\sigma}_i,
\]
where $(\innerprod{f}{\sigma_i})_{i \in I_N}$ are the $N$ largest
coefficients $\innerprod{f}{\sigma_i}$ in magnitude. As in the tight
frame case, this procedure does not always yield the {\em best}
$N$-term approximation, but, surprisingly, even with this rather crude
selection procedure, we can prove an (almost) optimally sparse
approximation rate. We speak of `almost' optimality due to the (negligible)
$\log$-factor in \eqref{eq:approx-rate-2d}. The following result shows
that our `new' compactly supported shearlets (see
Section~\ref{sec:constr-comp-supp}) deliver the same approximation
rate as \emph{band-limited} curvelets~\cite{CD04},
contourlets~\cite{DV05}, and shearlets~\cite{GL07}.

\begin{theorem}[\hspace*{-.25em}\cite{KL10}]
\label{theo:main}
Let $c > 0$, and let $\phi, \psi, \tilde{\psi} \in L^2(\R^2)$ be
compactly supported. Suppose that, in addition, for all $\xi =
(\xi_1,\xi_2) \in \R^2$, the shearlet $\psi$ satisfies
\begin{enumerate}[(i)]
\item $|\hat\psi(\xi)| \le C_1 \cdot \min\{1,|\xi_1|^{\alpha}\} \cdot
  \min\{1,|\xi_1|^{-\gamma}\} \cdot \min\{1,|\xi_2|^{-\gamma}\}$
  \hspace{2em} and
\item $\left|\frac{\partial}{\partial \xi_2}\hat \psi(\xi)\right| \le
  |h(\xi_1)| \cdot \left(1+\frac{|\xi_2|}{|\xi_1|}\right)^{-\gamma}$,
\end{enumerate}
where $\alpha > 5$, $\gamma \ge 4$, $h \in L^1(\R)$, and $C_1$ is a
constant, and suppose that the shearlet $\tilde{\psi}$ satisfies (i)
and (ii) with the roles of $\xi_1$ and $\xi_2$ reversed. Further,
suppose that $\SH(\phi,\psi,\tilde{\psi};c)$ forms a frame for
$L^2(\R^2)$.

Then, for any $\nu > 0$, the shearlet frame
$\SH(\phi,\psi,\tilde{\psi};c)$ provides (almost) optimally sparse
approximations of functions $f \in \cE^2_\nu(\R^2)$ in the sense that there exists
some $C > 0$ such that
\begin{equation}
\|f-f_N\|_2^2 \leq C\cdot N^{-2} \cdot {(\log{N})}^3 \qquad \text{as } N \rightarrow \infty,\label{eq:approx-rate-2d}
\end{equation}
where $f_N$ is the nonlinear N-term approximation obtained by choosing the N largest shearlet coefficients
of $f$.
\end{theorem}

Condition~(i) can be interpreted as both a condition ensuring (almost)
separable behavior as well as a moment condition along the horizontal
axis, hence enforcing directional selectivity. This condition ensures
that the support of shearlets in frequency domain is essentially of
the form indicated in Figure \ref{fig:2d-tiling}. Condition~(ii)
(together with (i)) is a weak version of a directional vanishing
moment condition\footnote{For the precise definition of directional
  vanishing moments, we refer to \cite{DV05}. }, which is crucial for
having fast decay of the shearlet coefficients when the corresponding
shearlet intersects the discontinuity curve. Conditions~(i) and (ii)
are rather mild conditions on the generators; in particular, shearlets
constructed by Theorem~\ref{thm:low-pass-construction}
and~\ref{thm:frame-for-R2}, with extra assumptions on the parameters
$K$ and $L$, will indeed satisfy (i) and (ii) in Theorem~\ref{theo:main}. To compare with the
optimality result for band-limited generators we wish to point out
that conditions~(i) and (ii) are obviously satisfied for band-limited
generators.

We remark that this kind of approximation result is not available for
shearlet systems coming directly from the shearlet group. One reason
for this being that these systems, as mentioned several times, do not
treat directions in a uniform way.

\section{Shearlets in 3D and Beyond}
\label{sec:shearl-high-dimens}
Shearlet theory has traditionally only dealt with representation
systems for two-dimensional data. In the recent paper \cite{DST09}
(and the accompanying paper \cite{DT10}) this was changed when Dahlke,
Steidl, and Teschke generalized the continuous shearlet transform (see
\cite{KL09,DKMSST08}) to higher dimensions. The shearlet transform on
$L^2(\R^n)$ by Dahlke, Steidl, and Teschke is associated with the
so-called shearlet group in \mbox{$\R \setminus \{ 0\} \times \R^{n-1}
\times \R^n$}, with a dilation matrix of the form
\[ A_a = \diag{a,\sign{a}\abs{a}^{1/n},\dots,\sign{a}\abs{a}^{1/n}},
\qquad a \in \R \setminus \{ 0\},\]
and with a shearing matrix with $n-1$ shear parameters $s=(s_1, \dots,
s_{n-1}) \in \R^{n-1}$ of the form
\[ S_s =
\begin{bmatrix}
  1 & s \\
  0 & I_{n-1}
\end{bmatrix},
\] where $I_n$ denotes the $n \times n$ identity matrix. This type of
shearing matrix gives rise to shearlets consisting of wedges of size
$a^{-1} \times a^{-1/n} \times \dots \times a^{-1/n}$ \emph{in
  frequency domain}, where $a^{-1} \gg a^{-1/n}$ for small $a>0$.
Hence, for small $a>0$, the spatial appearance is a surface-like
element of co-dimension one.

In the following section we will consider shearlet systems in
$L^2(\R^3)$ associated with a sightly different shearing matrix. More
importantly, we will consider \emph{pyramid-adapted} 3D shearlet
systems, since these systems treat directions in a uniform way as
opposed to the shearlet systems coming from the shearlet group; this
design, of course, parallels the idea behind cone-adapted 2D
shearlets. In \cite{GL10}, the continuous version of the
pyramid-adapted shearlet system was introduced, and it was shown that
the location and the local orientation of the boundary set of certain
three-dimensional solid regions can be precisely identified by this
continuous shearlet transform. The pyramid-adapted shearlet system can
easily be generalized to higher dimensions, but for brevity we only
consider the three-dimensional setup \gggk{and newly introduce it now
  in the discrete setting.}

\subsection{Pyramid-Adapted Shearlet Systems}
\label{sec:pyram-adapt-shearl}
We will scale according to \emph{paraboloidal scaling matrices}
$A_{2^j}$, $\tilde{A}_{2^j}$ or $\breve{A}_{2^j}$, $j \in \Z$, and
encode directionality by the \emph{shear matrices} $S_k$,
$\tilde{S}_k$, or $\breve{S}_k$, $k = (k_1,k_2) \in \Z^2$, defined by
\begin{alignat*}{6}
  A_{2^j} &=\begin{pmatrix}
    2^j & 0 & 0 \\ 0 & 2^{j/2} & 0 \\ 0 & 0 & 2^{j/2}
  \end{pmatrix}, &\quad&  \tilde{A}_{2^j}&=&\begin{pmatrix}
    2^{j/2} & 0 & 0 \\ 0 & 2^{j} & 0 \\ 0 & 0 & 2^{j/2}
  \end{pmatrix}, &\quad& \text{and}&\quad& \breve{A}_{2^j} &=&\begin{pmatrix}
    2^{j/2} & 0 & 0 \\ 0 & 2^{j/2} & 0 \\ 0 & 0 & 2^{j}
  \end{pmatrix},
\intertext{and}
S_k &=\begin{pmatrix} 1\; & k_1\; & k_2 \\ 0 & 1 & 0
    \\ 0 & 0 & 1
  \end{pmatrix}, &&  \tilde{S}_k &=&\begin{pmatrix} 1 & 0 & 0 \\ k_1\;
    & 1\; & k_2 \\ 0 & 0 & 1
  \end{pmatrix}, && \text{and} && \breve{S}_k\; &=&\begin{pmatrix} 1
    & 0 & 0 \\ 0 & 1 & 0 \\ k_1\; & k_2\; & 1
  \end{pmatrix},
\end{alignat*}
respectively. The translation lattices will be defined through the
following matrices: $M_c = \mathrm{diag}(c_1,c_2,c_2)$, $\tilde{M}_c =
\mathrm{diag}(c_2,c_1,c_2)$, and $\breve{M}_c =
\mathrm{diag}(c_2,c_2,c_1)$, where \jl{$c_1>0$ and $c_2>0$.}

\begin{figure}[ht]
\sidecaption[t]
\includegraphics[width=6cm]{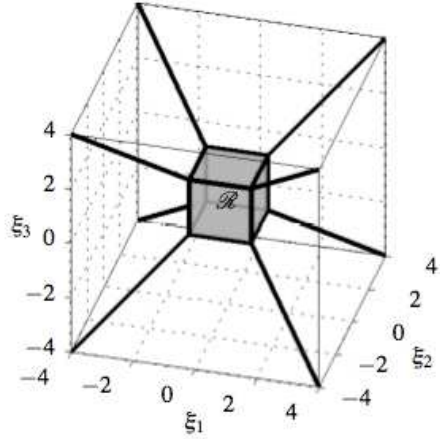}
      \caption{The partition of the frequency domain: The centered
        rectangle $\cR$. \jl{The arrangement of the six pyramids is
          indicated by the `diagonal' lines. See
          Figure~\ref{fig:pyramids} for a sketch of the pyramids.}}
\label{fig:partition}
\end{figure}

We next partition the frequency domain into the following six pyramids:
\begin{align*}
  \cP_\iota = \left\{ \begin{array}{rcl}
      \{(\xi_1,\xi_2,\xi_3) \in \R^3 : \xi_1 \ge 1,\, |\xi_2/\xi_1| \le 1,\, |\xi_3/\xi_1| \le 1\} & : & \iota = 1,\\
      \{(\xi_1,\xi_2,\xi_3) \in \R^3 : \xi_2 \ge 1,\, |\xi_1/\xi_2|
      \le 1,\,
      |\xi_3/\xi_2| \le 1\} & : & \iota = 2,\\
      \{(\xi_1,\xi_2,\xi_3) \in \R^3 : \xi_3 \ge 1,\, |\xi_1/\xi_3|
      \le 1,\,
      |\xi_2/\xi_3| \le 1\} & : & \iota = 3,\\
      \{(\xi_1,\xi_2,\xi_3) \in \R^3 : \xi_1 \le -1,\, |\xi_2/\xi_1| \le 1,\, |\xi_3/\xi_1| \le 1\} & : & \iota = 4,\\
      \{(\xi_1,\xi_2,\xi_3) \in \R^3 : \xi_2 \le -1,\, |\xi_1/\xi_2|
      \le
      1,\, |\xi_3/\xi_2| \le 1\} & : & \iota = 5,\\
      \{(\xi_1,\xi_2,\xi_3) \in \R^3 : \xi_3 \le -1,\, |\xi_1/\xi_3|
      \le 1,\, |\xi_2/\xi_3| \le 1\} & : & \iota = 6,
    \end{array}
  \right.
\end{align*}
and a centered rectangle
\[
\cR = \{(\xi_1,\xi_2,\xi_3) \in \R^3 : \norm[\infty]{(\xi_1,\xi_2, \xi_3)} < 1\}.
\]

\begin{figure}[ht]
 \vspace*{-.5em}
\centering
\subfloat[Pyramids $\cP_1$ and $\cP_4$ and the $\xi_1$ axis.]{%
\includegraphics[height=2.5cm]{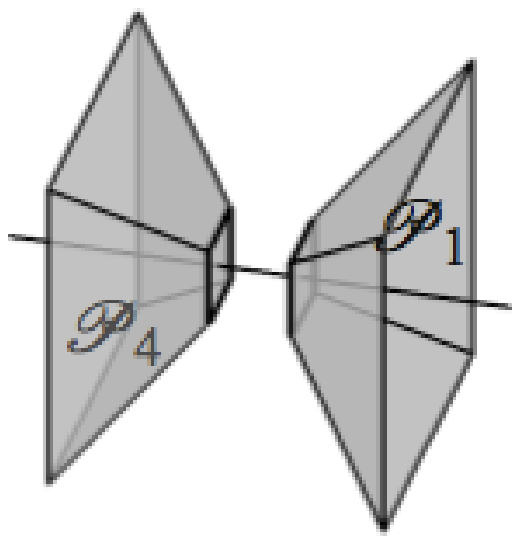}
\label{fig:pyramids14}
}\quad
\subfloat[Pyramids $\cP_2$ and $\cP_5$ and the $\xi_2$ axis.]{%
\includegraphics[height=2.5cm]{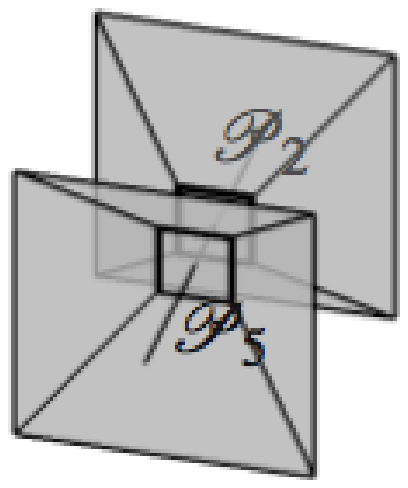}
\label{fig:pyramids25}
}\quad
\subfloat[Pyramids $\cP_3$ and $\cP_6$ and the $\xi_3$ axis.]{%
\includegraphics[height=2.5cm]{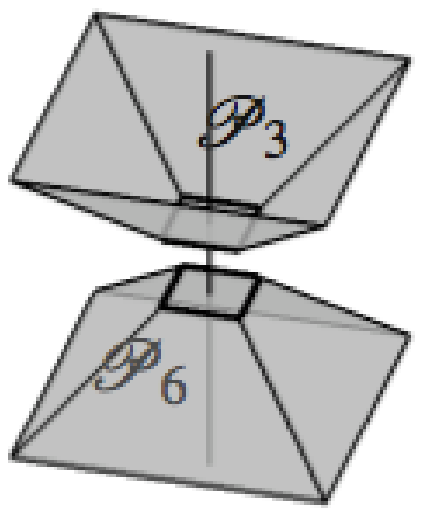}
\label{fig:pyramids136}
}
\caption{The partition of the
  frequency domain: The `top' of the six pyramids.}
\label{fig:pyramids}
\end{figure}

The partition is illustrated in Figures~\ref{fig:partition}
and~\ref{fig:pyramids}. This partition of the frequency space allows
us to restrict the range of the shear parameters. In the case of
`shearlet group' systems one must allow arbitrarily large shear
parameters, while the `pyramid-adapted' systems restrict the shear
parameters to $\itvcc{-\lceil2^{j/2}\rceil}{\lceil2^{j/2}\rceil}$. It
is exactly this fact that gives a more uniform treatment of the
directionality properties of the shearlet system.

These considerations are now made precise in the following definition.

\begin{definition}
  \label{def:discreteshearlets3d}
  For $c=(c_1,c_2) \in (\R_+)^2$, the \emph{pyramid-adapted 3D
    shearlet system} $\SH(\phi,\psi,\tilde{\psi},\breve{\psi};c)$
  generated by $\phi, \psi, \tilde{\psi}, \breve{\psi} \in L^2(\R^3)$
  is defined by
  \[
  \SH(\phi,\psi,\tilde{\psi},\breve{\psi};c) = \Phi(\phi;c_1) \cup
  \Psi(\psi;c) \cup \tilde{\Psi}(\tilde{\psi};c) \cup
  \breve{\Psi}(\breve{\psi};c),
  \]
  where
  \begin{align*}
    \Phi(\phi;c_1) &= \setprop{\phi_m = \phi(\cdot-m)}{ m \in c_1\Z^3}, \\
    \Psi(\psi;c) &= \setprop{\psi_{j,k,m} = \scfac {\psi}({S}_{k}
    {A}_{2^j}\cdot-m) }{ j \ge 0, |k| \le \ceilsmall{2^{j/2}}, m \in
    M_c \Z^3 }, \\
    \tilde{\Psi}(\tilde{\psi};c) &= \{\tilde{\psi}_{j,k,m} = \scfac
    \tilde{\psi}(\tilde{S}_{k} \tilde{A}_{2^j}\cdot-m) : j \ge 0, |k| \le
    \lceil 2^{j/2} \rceil, m \in \tilde{M}_c \Z^3 \},\\
    \intertext{and} \breve{\Psi}(\breve{\psi};c) &=
    \{\breve{\psi}_{j,k,m} = \scfac \breve{\psi}(\breve{S}_{k}
    \breve{A}_{2^j}\cdot-m) : j \ge 0, |k| \le \lceil 2^{j/2} \rceil,
    m \in \breve{M}_c \Z^3 \},
  \end{align*}
  where $j \in \N_0$ and $k \in \Z^2$. 
  Here we have used the vector
  notation $\abs{k} \le K$ for $k = (k_1,k_2)$ and $K>0$ to denote
  $\abs{k_1} \le K$ \emph{and} $\abs{k_2} \le K$.
\end{definition}

The construction of pyramid-adapted shearlet systems
$\SH(\phi,\psi,\tilde{\psi},\breve{\psi};c)$ runs along the lines of
the construction of cone-adapted shearlet systems in $L^2(\R^2)$
described in Section~\ref{sec:constr-comp-supp}. For a detailed
description, we refer to \cite{KLL10}.

We remark that the shearlets in {\em spatial domain} are of size
$2^{-j/2}$ times $2^{-j/2}$ times $2^{-j}$ which shows that the
shearlet elements will become `plate-like' as $j \to \infty$. One
could also use the scaling matrix $A_{2^j}=\diag{2^j,2^j,2^{j/2}}$ with
similar changes for $\tilde A_{2^j}$ and $ \breve A_{2^j}$. This would
lead to `needle-like' shearlet elements instead of the `plate-like'
elements considered in this paper, but we will not pursue this further
here, and simply refer to \cite{KLL10}. More generally, it is
possible to even consider non-paraboloidal scaling matrices of the form
$A_{j}=\diag{2^j,2^{\alpha j},2^{\beta j}}$ for $0 < \alpha, \beta \le
1$. One drawback of allowing such general scaling matrices is the lack
of fast algorithms for non-dyadic multiscale systems. On the other
hand, the parameters $\alpha$ and $\beta$ allow us to precisely shape
the shearlet elements, ranging from very plate-like to very
needle-like, according to the application at hand, \ie choosing the
shearlet-shape that is the best `fit' for the geometric characteristics
of the considered data.

\subsection{Sparse Approximations of 3D Data}
\label{sec:sparse-appr-3d}

We now consider approximations of three-dimensional cartoon-like
images using shearlets introduced in the previous section. The
three-dimensional cartoon-like images $\cE^2_\nu(\R^3)$ will be
piecewise $C^2$
functions with discontinuities on a closed $C^2$ \emph{surface} whose
principal curvatures are bounded by $\nu$. In \cite{KLL10} it was shown that the optimal
approximation rate for such 3D cartoon-like image models $f \in
\cE^2_\nu(\R^3)$ which can be achieved for almost any representation
system (under polynomial depth search selection procedure of the
approximating coefficients) is
\[
\norm[2]{f-f_N}^2 \le C \cdot N^{-1} \quad \text{as } N \to \infty,
\]
where $f_N$ is the best $N$-term approximation of $f$. The following
result shows that compactly supported pyramid-adapted shearlets do
(almost)  deliver this approximation rate.

 \begin{theorem}[\hspace*{-.24em}\cite{KLL10}]
\label{thm:3d-opt-sparse}
Let $c \in (\R_+)^2$, and let $\phi, \psi, \tilde{\psi}, \breve{\psi}
\in L^2(\R^3)$ be compactly supported. Suppose that, for all $\xi =
(\xi_1,\xi_2,\xi_3) \in \R^3$, the function $\psi$ satisfies:
   \begin{enumerate}[(i)]
   \item $|\hat\psi(\xi)| \le C_1 \cdot \min\{1,|\xi_1|^{\alpha}\} \cdot
     \min\{1,|\xi_1|^{-\gamma}\} \cdot \min\{1,|\xi_2|^{-\gamma}\} \cdot
     \min\{1,|\xi_3|^{-\gamma}\}$,
   \item $\left|\frac{\partial}{\partial \xi_i}\hat \psi(\xi)\right|
     \le |h(\xi_1)| \cdot
     \left(1+\frac{|\xi_2|}{|\xi_1|}\right)^{-\gamma}
     \left(1+\frac{|\xi_3|}{|\xi_1|}\right)^{-\gamma}, \qquad i=2,3$,
   \end{enumerate}
   where $\alpha > 8$, $\gamma \ge 4$, $t \mapsto th(t) \in L^1(\R)$, and $C_1$ a
   constant, and suppose that $\tilde{\psi}$ and $\breve{\psi}$
   satisfy analogous conditions with the obvious change of
   coordinates. Further, suppose that the shearlet system
   $\SH(\phi,\psi,\tilde{\psi},\breve{\psi};c)$ forms a frame for
   $L^2(\R^3)$.

   Then, for any $\nu > 0$, the shearlet frame
   $\SH(\phi,\psi,\tilde{\psi},\breve{\psi};c)$ provides (almost) optimally sparse
   approximations of functions $f \in \cE^2_\nu(\R^3)$ in the sense that there
   exists some $C > 0$ such that
   \begin{equation}\label{eq:3d-approx-rate}
     \|f-f_N\|_2^2 \leq C\cdot N^{-1} \cdot {(\log{N})}^2 \qquad \text{as } N \rightarrow \infty.
\end{equation}
\end{theorem}

 In the following we will give a sketch of the proof of
 Theorem~\ref{thm:3d-opt-sparse} and, in particular, give a
 heuristic argument (inspired by a similar one for 2D curvelets in
 \cite{CD04}) to explain the exponent $N^{-1}$ in
 (\ref{eq:3d-approx-rate}).

\begin{proof}[Theorem~\ref{thm:3d-opt-sparse}, Sketch]
  Let $f \in \cE^2_{\nu}(\R^3)$ be a 3D cartoon-like image. The main
  concern is to derive appropriate estimates for the shearlet
  coefficients $\innerprod{f}{\psi_{j,k,m}}$. We first observe that we
  can assume the scaling index $j$ to be sufficiently large, since $f$
  as well as all shearlet elements are compactly supported and since a
  finite number does not contribute to the asymptotic estimate we are
  aiming for. In particular, this implies that we do not need to take
  frame elements from the `scaling' system $\Phi(\phi;c_1)$ into
  account. Also, we are allowed to restrict our analysis to shearlets
  $\psi_{j,k,m}$, since the frame elements $\tilde{\psi}_{j,k,m}$ and
  $\breve{\psi}_{j,k,m}$ can be handled in a similar way.

  Letting $|\theta(f)|_n$ denote the $n$th largest shearlet
  coefficient $\innerprod{f}{\psi_{j,k,m}}$ in absolute value and
  using the frame property of
  $\SH(\phi,\psi,\tilde{\psi},\breve{\psi};c)$, we conclude that
  \[
  \|f-f_N\|_2^2 \leq \frac{1}{A}\sum_{n>N} |\theta(f)|_n^2,
  \]
  for any positive integer $N$, where $A$ denotes the lower frame
  bound of the shearlet frame
  $\SH(\phi,\psi,\tilde{\psi},\breve{\psi};c)$. Thus, for completing
  the proof, it therefore suffices to show that
  \begin{equation}\label{eq:upper}
    \sum_{n>N} |\theta(f)|_n^2 \leq C \cdot N^{-1} \cdot  {(\log{N})}^2 \qquad \text{as} \,\, N \rightarrow \infty.
  \end{equation}

  For the following heuristic argument, we need to make some
  simplifications. We will assume to have a shearlet of the form
  $\psi(x)=\eta(x_1)\varphi(x_2)\varphi(x_3)$, where $\eta$ is a
  wavelet and $\varphi$ a bump (or a scaling) function. Note that the
  wavelet `points' in the short direction of the plate-like shearlet.
  We now consider three cases of coefficients
  $\innerprod{f}{\psi_{j,k,m}}$:
  \begin{enumerate}[(a)]
  \item Shearlets \gggk{$\psi_{j,k,m}$} whose support does not overlap
    with the boundary $\partial B$.
  \item Shearlets \gggk{$\psi_{j,k,m}$} whose support overlaps with $\partial B$
    and is nearly tangent.
  \item Shearlets \gggk{$\psi_{j,k,m}$} whose support overlaps with
    $\partial B$, but not tangentially.
  \end{enumerate}

\begin{figure}[ht]
\centering
\subfloat[Sketch of shearlets whose support does not overlap with
    $\partial B$.]{%
\parbox[t]{0.29\textwidth}{\centering\includegraphics[height=3cm]{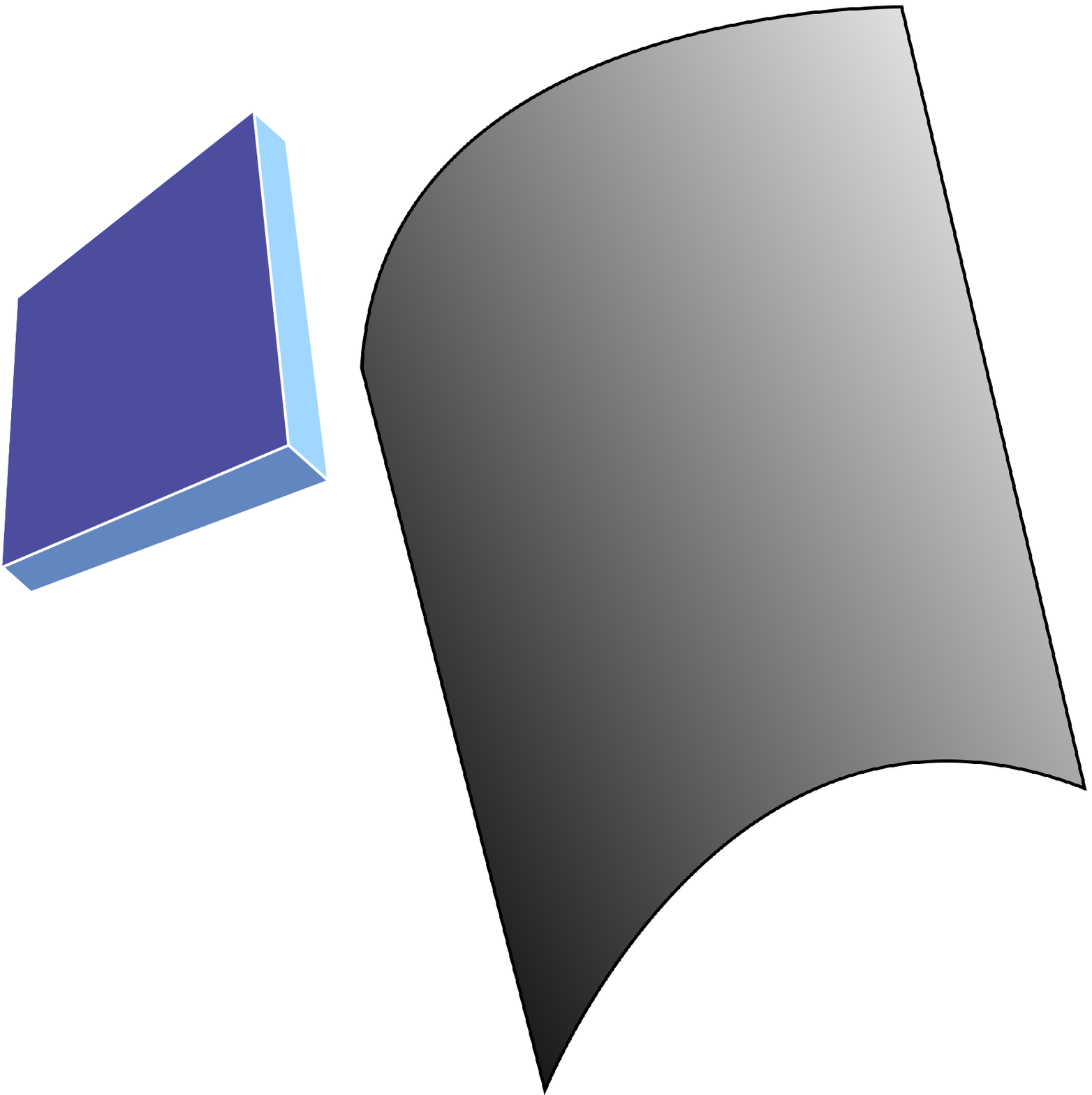}}
\label{fig:case1}
}\qquad
\subfloat[Sketch of shearlets whose support overlaps with
    $\partial B$ and is nearly tangent.]{%
\parbox[t]{0.29\textwidth}{\centering\includegraphics[height=3cm]{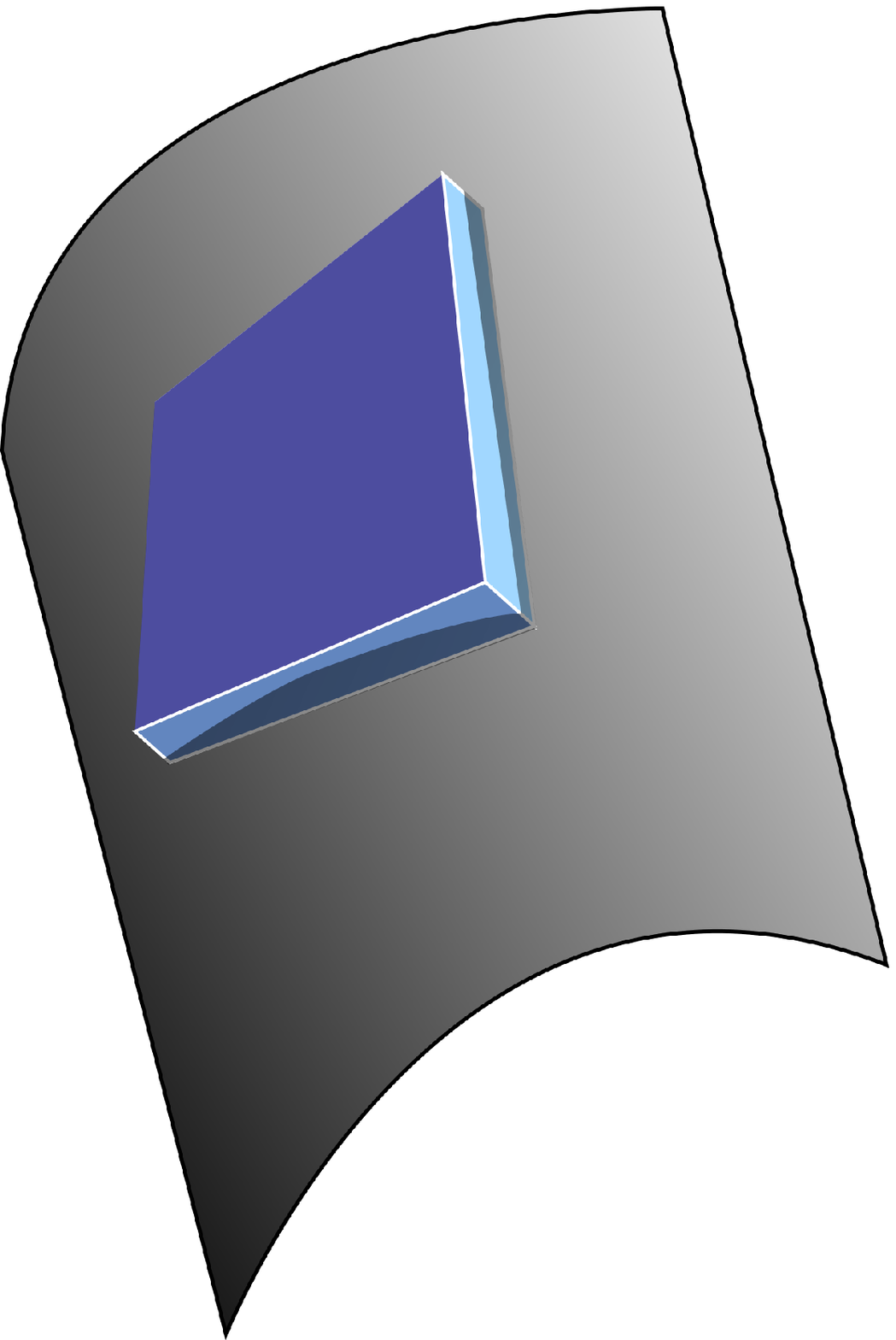}}
\label{fig:case2}
}\qquad
\subfloat[Sketch of shearlets whose support overlaps with
    $\partial B$, but not tangentially.]{%
\parbox[t]{0.29\textwidth}{\centering\includegraphics[height=3cm]{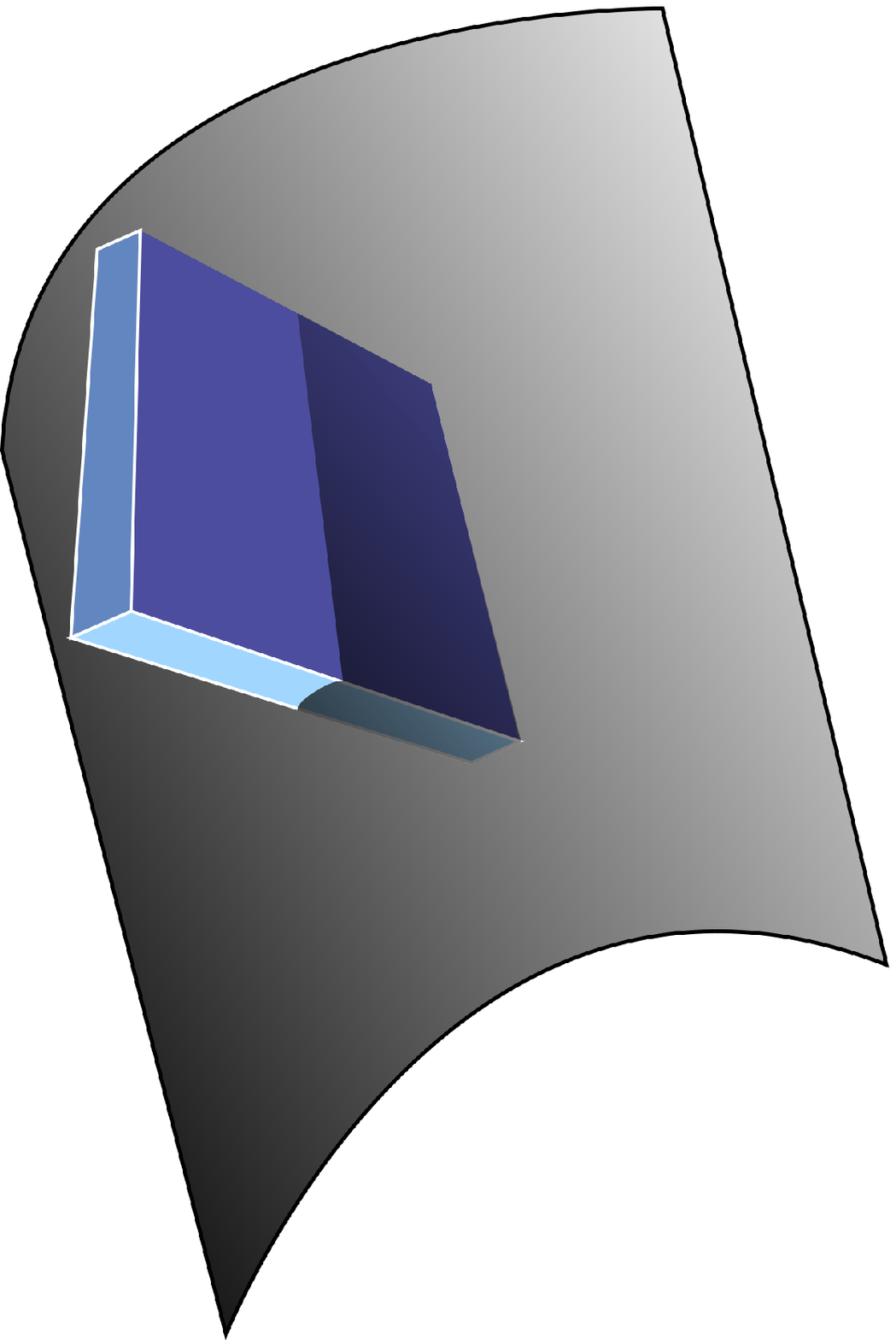}}
\label{fig:case3}
}
\caption{The three types of shearlet \gggk{$\psi_{j,k,m}$} and boundary $\partial B$
  interactions considered in the heuristic argument (explaining the
  approximation rate $N^{-1}$). Note that only a section of $\partial
  B$ is shown.}
\label{fig:three-cases}
\end{figure}
As we argue in the following, only coefficients from case~(b) will be
significant. Case~(b) is -- loosely speaking -- the situation in which
the wavelet $\eta$ \gggk{breaches}, in an almost normal direction,
through the discontinuity surface; as is well known from wavelet
theory, 1D wavelets efficiently handle such a `jump' discontinuity.

\emph{Case (a).} Since $f$ is $C^2$ smooth away from $\partial B$, the
coefficients $|\innerprod{f}{\psi_{j,k,m}}|$ will be sufficiently
small owing to the wavelet $\eta$ (and the fast decay of wavelet
coefficients of smooth functions).

\emph{Case (b).} At scale $j>0$, there are at most $O(2^{j})$
coefficients, since the plate-like elements are of size $2^{-j/2}$ times
$2^{-j/2}$ (and `thickness' $2^{-j}$). By assumptions on $f$ and the
support size of $\psi_{j,k,m}$, we obtain the estimate
\[ |\innerprod{f}{\psi_{j,k,m}}| \le \norm[\infty]{f}
  \norm[1]{\psi_{j,k,m}} \le C_1 \, (2^{-2j})^{1/2}
  \norm[2]{\psi_{j,k,m}}^{1/2} \le C_2 \cdot 2^{-j} \]
for some constants $C_1, C_2>0$. In other words, we have $O(2^{j})$
coefficients bounded by $C_2 \cdot 2^{-j}$. Assuming the case (a) and
(c) coefficients are negligible, the $n$th largest coefficient
$|\theta(f)|_n$ is then bounded by
\[ |\theta(f)|_n \le C \cdot n^{-1}.\]
Therefore,
\[
    \sum_{n>N} |\theta(f)|_n^2 \leq \sum_{n>N}C \cdot n^{-2} \le
    C \cdot \int_N^\infty x^{-2} dx \le C \cdot N^{-1}
 \]
 and we arrive at (\ref{eq:upper}), but without the $\log$-factor.
 This in turn shows (\ref{eq:3d-approx-rate}), at least heuristically,
 and still without the $\log$-factor.

\emph{Case (c).} Finally, when the shearlets are sheared away from the
tangent position in case~(b), they will again be small. This is due to
the vanishing moment conditions in condition (i) and (ii).
\qed
\end{proof}

Clearly, Theorem~\ref{thm:3d-opt-sparse} is an `obvious'
three-dimensional version of Theorem~\ref{theo:main}. However, as
opposed to the two-dimensional setting, anisotropic structures in
three-di\-men\-si\-onal data comprise of \emph{two} morphologically
different types of structure, namely surfaces \emph{and} curves. It
would therefore be desirable to allow our 3D image class to also
contain cartoon-like images with \emph{curve} singularities. On the
other hand, the pyramid-adapted shearlets introduced in
Section~\ref{sec:pyram-adapt-shearl} are plate-like and thus, a
priori, not optimal for capturing such one-dimensional singularities.
Surprisingly, these plate-like shearlet systems still deliver the
optimal rate $N^{-1}$ for three-dimensional cartoon-like images
$\cE^2_{\nu,L}(\R^3)$, where $L$ indicates that we allow our
discontinuity surface $\partial B$ to be \emph{piecewise} $C^2$
smooth; $L \in \N$ is the \gggk{maximal} number of $C^2$ pieces and
$\nu > 0$ is an upper estimate for the principal curvatures on each
piece. In other words, for any $\nu > 0$ and $L \in \N$, the shearlet
frame $\SH(\phi,\psi,\tilde{\psi},\breve{\psi};c)$ provides (almost)
optimally sparse approximations of functions $f \in
\cE^2_{\nu,L}(\R^3)$ in the sense that there exists some $C > 0$ such
that
\begin{equation}\label{eq:3d-approx-rate2}
  \|f-f_N\|_2^2 \leq C\cdot N^{-1} \cdot {(\log{N})}^2 \qquad \text{as } N \rightarrow \infty.
\end{equation}

The conditions on the shearlets $\psi, \tilde{\psi}, \breve{\psi}$
are similar to these in Theorem~\ref{thm:3d-opt-sparse}, but more
technical, and we refer to \cite{KLL10} for the precise statements and
definitions as well as the proof of the optimal approximation error rate.
Here we simply remark that there exist numerous examples of shearlets $\psi, \tilde{\psi}$,
and $\breve{\psi}$ satisfying these conditions, which lead to \eqref{eq:3d-approx-rate2}; one large class of examples
are separable generators $\psi, \tilde{\psi}, \breve{\psi} \in
L^2(\R^3)$, i.e.,
\[\psi(x)=\eta(x_1) \varphi(x_2)\varphi(x_3), \quad
\tilde\psi(x)=\varphi(x_1) \eta(x_2)\varphi(x_3), \quad
\breve\psi(x)=\varphi(x_1) \varphi(x_2)\eta(x_3),\] where $\eta,
\varphi \in L^2(\R)$ are compactly supported functions satisfying:
   \begin{enumerate}[(i)]
   \item $|\hat\eta(\omega)| \le C_1 \cdot \min\{1,|\omega|^{\alpha}\}
     \cdot \min\{1,|\omega|^{-\gamma}\}$, \vspace{0.2em}
   \item $\abs{\bigl(\frac{\partial}{\partial \omega}\bigr)^{\ell}\hat
       \varphi(\omega)} \leq C_2 \cdot \min\{1,|\omega|^{-\gamma}\}
     \quad$ for $\ell=0,1$,
   \end{enumerate}
   for $\omega \in \R$, where $\alpha > 8$, $\gamma \ge 4$, and $C_1,
   C_2$ are constants.

\section{Conclusions}
\label{sec:conclusion}
Designing a directional representation system that efficiently handles
data with anisotropic features is quite challenging since it needs to
satisfy a long list of desired properties: it should have a simple
mathematical structure, it should provide optimally sparse
approximations of certain image classes, it should allow compactly
supported generators, it should be associated with fast decomposition
algorithms, and it should provide a unified treatment of the continuum
and digital realm.

In this paper, we argue that shearlets meet all these challenges, and
are, therefore, one of the most satisfying directional systems. To be
more precise, let us briefly review our findings for 2D and 3D data:
\begin{itemize}
\item \emph{2D Data.} In Section~\ref{sec:shearl-two-dimens}, we
  constructed 2D shearlet systems that efficiently capture anisotropic
  features and satisfy all the above requirements.
\item \emph{3D Data.} In 3D, as opposed to 2D, we face the difficulty
  that there might exist two geometrically different anisotropic
  features; 1D and 2D singularities. The main difficulty in extending
  shearlet systems from the 2D to 3D setting lies, therefore, in
  introducing a system that is able to represent both these
  geometrically different structures efficiently. As shown in
  Section~\ref{sec:shearl-high-dimens}, a class of plate-like
  shearlets is able to meet these requirements. In other words: the
  extension from 2D shearlets to 3D shearlets has been successful in
  terms of preserving the desirable properties, \eg optimally sparse
  approximations. It does therefore seem that an extension to 4D or
  even higher dimensions is, if not straightforward then, at the very
  least, feasible. In particular, the step to 4D now `only' requires
  the efficient handling of yet `another' type of anisotropic feature.
\end{itemize}


\begin{acknowledgement}
  The first and third author acknowledge support from DFG Grant SPP-1324, KU 1446/13.
  The first author also acknowledges support from DFG Grant KU 1446/14.
\end{acknowledgement}

\end{document}